\documentclass[12pt,twoside,leqno]{amsart}

\usepackage{amssymb}
\usepackage{amsbsy}
\usepackage{amsmath}
\usepackage[english]{babel}
\usepackage{multicol}
\usepackage{subfig}
\usepackage{times}
\usepackage[all]{xy}
\usepackage{comment}
\usepackage{cyrillic} 
\usepackage[T1]{fontenc} 

\newtheorem{proposition}{Proposition}[section]
\newtheorem{corollary}[proposition]{Corollary}
\newtheorem{theorem}[proposition]{Theorem}
\newtheorem{recall}[proposition]{Recall}
\theoremstyle{definition}
\newtheorem{definition}[proposition]{Definition}
\newtheorem{lemma}[proposition]{Lemma}

\theoremstyle{remark}
\newtheorem{remark}[proposition]{Remark}

\newcommand{\cyrrm}{\fontencoding{OT2}\selectfont\textcyrup}

\DeclareMathAlphabet{\mathbbold}{U}{bbold}{m}{n}

\title{Borel-Moore motivic homology and weight structure on mixed motives}

\author[]{Jin Fangzhou}
\date{\number\day-\number\month-\number\year}

\begin{document}

\maketitle

\begin{abstract}
By defining and studying functorial properties of the Borel-Moore motivic homology, we identify the heart of Bondarko-H\'ebert's weight structure on Beilinson motives with Corti-Hanamura's category of Chow motives over a base, therefore answering a question of Bondarko.

\end{abstract}

\tableofcontents

\noindent

\section{Introduction}


Bondarko introduced in \cite{Bon1} the notion of \emph{weight structure} on triangulated categories as a counterpart of $t$-structure (\cite{BBD}), with the aim of applying it to the theory of mixed motives. A weight structure defines a collection of non-positive and non-negative objects, such that every object can be decomposed, via a distinguished triangle, into a non-positive part and a non-negative part, in a non-unique way. The weight-$0$ part, also called the heart, is formed by objects which are both non-positive and non-negative; it carries so much information that a weight structure can be recovered from a partial knowledge on its heart: to be precise, any additive subset of non-positive objects is contained in the heart of a unique weight structure on the subcategory it generates (\cite[Proposition 5.2.2]{Bon1}). As an application of this abstract notion, Bondarko showed the existence of a weight structure on Voevodsky's motives over a perfect field (\cite{Voe}), which provides a motivic analogue of Deligne's weight filtration and weight spectral sequence (\cite{Del}). He also showed (\cite[6.5 and 6.6]{Bon1} ) that the heart is the category of Chow motives over a field (\cite[4.1.3]{And} ).

Based on Ayoub's work on cross functors in \cite{Ayo}, Cisinski-D\'eglise generalized Voevodsky's motivic complexes to arbitrary base scheme in \cite{CD}, by constructing the category of Beilinson motives:
for any scheme $S$, there is a $\mathbb{Q}$-linear triangulated category $DM_{{\cyrrm{B}}}(S)$ in which the \emph{Grothendieck six functors formalism} is satisfied (\emph{ibid.} A.5.1., and we will recall part of the full formalism in Recall~\ref{recall_DM}). The objects of geometric nature, namely those which satisfy some finiteness condition, are called constructible motives, and form a full subcategory $DM_{{\cyrrm{B}},c}(S)$ (\emph{ibid.} C.2.). In the case of Beilinson motives, H\'ebert (\cite{Heb}) and independently Bondarko (\cite{Bon2}) showed, by different methods, the existence of weight structures on both $DM_{{\cyrrm{B}}}(S)$ and $DM_{{\cyrrm{B}},c}(S)$.

The aim of the present article is to give a detailed description of the heart in $DM_{{\cyrrm{B}},c}(S)$, namely answer the conjecture in \cite[Remark 2.1.2.4]{Bon2} : if $S$ is a quasi-projective scheme over a perfect field $k$, the heart of this weight-structure, denoted by $Chow(S)$, is equivalent to the category of Chow motives over a base $CHM(S)$, defined in \cite[Definition 2.8]{CH}. Recall that the latter is constructed in a similar way as Chow motives over a field, where the composition of cycles is defined in the following way: let $X$ (respectively $Y$, $Z$) be a scheme endowed with a projective morphism $f:X\to S$ such that $X$ is smooth of relative dimension $d_X$ (respectively $d_Y$, $d_Z$) over $k$. Then the following map between Chow groups defines composition between correspondences:
\begin{align*}
       CH_{d_Y+n}(X\times_SY) \times CH_{d_Z+m}(Y\times_SZ) &\to CH_{d_Z+m+n}(X\times_SZ)\\
        (\alpha,\beta)&\mapsto p^{XYZ}_{XZ*}\ \delta_Y^!(\alpha\times\beta),
\end{align*}
where $p^{XYZ}_{XZ}$ is the projection $X\times_SY\times_SZ\to X\times_SZ$, $\alpha\times\beta$ is the exterior product (\cite[1.10]{Ful} ), and $\delta_Y^!$ is the refined Gysin morphism (\emph{ibid.} 6.2.) associated to the cartesian diagram
$$
  \xymatrix{
    X\times_SY\times_SZ \ar[r] \ar[d] & (X\times_SY)\times_k(Y\times_SZ) \ar[d]\\
    Y \ar[r]^-{\delta_Y} & Y\times_kY,
  }
$$
where the diagonal embedding $\delta_Y$ is a regular embedding (\emph{ibid.} B.7.) since $Y$ is smooth over $k$.
There is a canonical map from $CHM(S)$ to the heart, which sends the motive of $X$ to the element $M^{BM}(X/S)=f_!\mathbbold{1}_X$ (where $\mathbbold{1}_X$ is the unit for the tensor structure), called the \emph{Borel-Moore motive} of $X$ relative to $S$. A natural idea is to show that this map is indeed a functor, which induces an equivalence of categories. To do so, we need to
\begin{enumerate}
\item Identify the group of morphisms between Borel-Moore motives;
\item Describe the composition law between morphisms.
\end{enumerate}

For the first point, using a base change formula, we are reduced to compute the following group
$$
H^{BM}_{p,q}(X/S)=Hom_{DM_{{\cyrrm{B}},c}(S)}(M^{BM}(X/S),\mathbbold{1}_S(-q)[-p])
$$
(see also \cite[Lemma 1.1.4]{Bon2}), which we call the \emph{Borel-Moore motivic homology} (or simply \emph{Borel-Moore homology}) of $X$ relative to $S$. 
Classically, for $S=Spec\ \mathbb{C}$, Borel-Moore homology, also known as homology with compact support, is a homology theory for locally compact topological spaces introduced 
in \cite{BM}. A motivic analogue is defined by Levine (\cite[7.8.5]{Lev2} ) using the notion of motives with compact support, and a similar definition is given in \cite[Definition 16.20]{MVW} . Our treatment here uses the powerful tool of six functors formalism, which provides similar and even stronger functorialities compared with the topological case, such as the following ones, called "basic functorialities":
\begin{itemize}
\item For any proper morphism $f:X\to Y$, there is a push-forward map $$f_*:H^{BM}_{p,q}(X/S)\overset{}{\to}H^{BM}_{p,q}(Y/S);$$
\item For any smooth morphism $f:X\to Y$ of relative dimension $d$, there is a pull-back map $$f^*:H^{BM}_{p,q}(Y/S)\overset{}{\to}H^{BM}_{p+2d,q+d}(X/S);$$
\item For any closed immersion $i:Z\to X$ with complement open immersion $j:U\to X$, there is a long exact sequence, called \emph{localization sequence}
$$
\cdots\to
H^{BM}_{p,q}(Z/S)\overset{i_*}{\to}
H^{BM}_{p,q}(X/S)\overset{j^*}{\to}
H^{BM}_{p,q}(U/S)\overset{\partial_{X,Z}}{\to}
H^{BM}_{p-1,q}(Z/S)
\to\cdots.
$$
\end{itemize}
Computations of Borel-Moore homology date back to Voevodsky's identification of motivic cohomology for smooth schemes with Chow groups in \cite{Voe2}. 
The key in his proof is the use of the \emph{coniveau spectral sequence}, introduced for the first time by Bloch and Ogus for \'etale cohomology (\cite{BO}) after Grothendieck's insight on the coniveau filtration (\cite{Gro}). A big amount of work has been devoted to this kind of spectral sequence since their work: for example \cite{CTHK} computed the $E_2$ terms for more general cohomologies; Rost introduced the notion of \emph{cycle modules} in \cite{Ros} in an attempt to deal with a generalized divisor class map from the prototype of Milnor's K-theory; D\'eglise established the link between Voevodsky's motives and Rost's cycle modules in \cite{Deg7} and \cite{Deg4}, and developed the corresponding motivic theory in \cite{Deg1}. The $E_1$-term of the coniveau spectral sequence, also known as \emph{Cousin complexes}, come from filtration by codimension of support: they are sums of "generic" cohomology groups corresponding to points of $X$ of a given codimension, where each point is seen as the generic point of its Zariski closure. 
The coniveau spectral sequence converges to the cohomology of the entire variety, which in our case allows to compute certain motivic cohomology groups. Another crucial ingredient in Voevodsky's proof is the cancellation of certain motivic cohomology groups:
\begin{lemma}
\label{lemma_voe}
(\cite[Lemma 3.2]{SV})
Let $k$ be a field and $X$ be a smooth $k$-scheme, then we have
$$
H^{p,q}(X)=0
$$
whenever $p>q+dim(X)$ or $q<0$.
\footnote{We will define motivic cohomology groups in Definition~\ref{def_mot_coh}.}
\end{lemma}
As a consequence, the index of non-zero $E_1$-terms is constrained in a quarter-plane of $\mathbb{Z}^2$; identifying the differential map with the divisor class map (see for example \cite[Proposition 1.16]{Deg1} ), we obtain the following
\begin{lemma}
(\cite[Proposition 2.1.4]{Voe})
Let $k$ be a perfect field and $X$ be a smooth $k$-scheme, then there is an isomorphism
$$
H^{2n,n}(X)\simeq CH^n(X)\otimes \mathbb{Q}.
$$
\end{lemma}
As a matter of fact, the localization property of Chow groups suggests that they behave more like a Borel-Moore theory. For smooth schemes, motivic cohomology agrees with Borel-Moore homology up to a suitable change of indexes (Lemma~\ref{ident_BM_mot}); this agreement fails for singular schemes, in which case it is the Borel-Moore homology that gives the correct Chow groups, as shown in \cite[Proposition 19.18]{MVW}. We will give a proof in our setting by constructing a spectral sequence similar to the coniveau spectral sequence which, as we will see, is essentially a consequence of the localization sequence; using some results on elementary operations (Lemma~\ref{elem_BM}), we deduce the following
\begin{proposition}
(see Corollary~\ref{BM=Chow})
Let $k$ be a perfect field and $X$ be a separated $k$-scheme of finite type, then there is an isomorphism
$$
H^{BM}_{2n,n}(X)\simeq CH_n(X)\mathbb{Q}.
$$
\end{proposition}





For the second point, seeing the composition law on Chow motives, it is important to construct a motivic version of Fulton's refined Gysin morphism, and then compare functorialities on Chow groups and those on motives. There are indeed two ways to define refined Gysin morphisms on motives: on the one hand, using the absolute purity isomorphism for closed immersion between regular schemes, we show that there is a natural way to construct Gysin morphisms, which are generalized to refined Gysin morphisms by base change; on the other hand, we can repeat the construction in \cite{Ful} by constructing a motivic specialization map and then composing with suitable push-forward and pull-back maps. By a non-trivial lemma, we show that the two possible definitions agree (Corollary~\ref{refined_gysin_compat}). This two-sided point of view is very advantageous: the first definition is so natural that it can be easily compared with the composition law with a formal argument (Proposition~\ref{composition_law}); the second one can be compared with Fulton's construction, if we know how to compare push-forward, pull back and specialization maps. This can be done by looking directly at the niveau spectral sequence. To conclude, the main result that we obtain is the following:
\begin{theorem}
(see Theorem~\ref{main_thm})
Let $S$ be a quasi-projective scheme over a perfect field $k$. Then the two categories $CHM(S)$  and $Chow(S)$ are equivalent.
\end{theorem}





\medskip
\noindent
The article is organized as follows:

In Chapter~\ref{chapter2} we study several types of functorialities of the Borel-Moore motive. Apart from the "basic functorialities" listed above, purity isomorphisms provide much important information: we use them to define Gysin morphisms and therefore deduce refined Gysin morphisms by base change. By a similar lemma as \cite[Proposition 2.6.5]{Deg4}, we show that this definition agrees with the construction in \cite[6.2]{Ful}. At the end of the chapter, we show that refined Gysin morphisms we have defined appear naturally when it comes to the composition of morphisms between Borel-Moore motives.


In Chapter~\ref{Chapter3} we begin with the computation of the Borel-Moore homology. We construct the so-called niveau spectral sequence, similar to the coniveau one, for Borel-Moore homology; as a consequence, we deduce that Borel-Moore homology groups are isomorphic to Chow groups (Corollary~\ref{BM=Chow}); then we show the compatibility of this isomorphism with respect to push-forward, pull-back and refined Gysin maps. At last we conclude by giving the identification between Chow motives over a base and the heart of Bondarko-H\'ebert's weight structure.

In the last chapter we prove some lemmas used in previous chapters. Section~\ref{proof_elem_BM} is devoted to the proof of Lemma~\ref{elem_BM}, which computes some elementary operations on Borel-Moore homology. Such computations are already known in \cite{Deg2} for motivic cohomology, and we adapt those results to our case by identifying motivic cohomology with Borel-Moore homology. In Section~\ref{proof_BC_proper_smooth} we prove the remaining lemmas via a further study of purity isomorphisms.

\medskip\noindent\textbf{Acknowledgments.}\ \ I would like to express my gratitude to my thesis advisor, Fr\'ed\'eric D\'eglise, who kindly shared the basic ideas of this article and helped me greatly during its development. I would also like to thank Mikhail Bondarko for very helpful discussion and advice on a preprint version.

\section{Borel-Moore motives and their functorialities}
\label{chapter2}

In this chapter we study functorial properties of Borel-Moore motives. Note that the proof of Lemma~\ref{BC_proper_smooth}, Lemma~\ref{compat_purity}, Lemma~\ref{compat_gysin_func} and Lemma~\ref{lemma_monoidal_BC} will be postponed to Section~\ref{proof_BC_proper_smooth}.

\smallskip
\noindent
\textbf{Convention.}\ If $(F,G)$ is a pair of adjoint functors, we denote by $ad_{(F,G)}:1\to G\circ F$ and $ad'_{(F,G)}:F\circ G\to 1$ the canonical counit and unit natural transformations.

\subsection{The category of Beilinson motives}\

Following \cite{CD}, the category of constructible Beilinson motives, is a triangulated category $DM_{{\cyrrm{B}},c}$ fibered over the category of noetherian schemes of finite Krull dimension, 
which satisfies the 6 functors formalism. Recall some of its properties, which we will need later:
\begin{recall}
\label{recall_DM}
In the category $DM_{{\cyrrm{B}},c}$, the following conditions hold:
\begin{enumerate}
\item There is symmetric closed monoidal structure on $DM_{{\cyrrm{B}},c}$, denoted by $\otimes$. For any $X\in\mathcal{S}$, the unit in $DM_{{\cyrrm{B}},c}(X)$ is denoted by $\mathbbold{1}_X$.
\item For any morphism $f:Y\to X$, there is a functor $f^*:DM_{{\cyrrm{B}},c}(X)\to DM_{{\cyrrm{B}},c}(Y)$, which has a right adjoint $f_*:DM_{{\cyrrm{B}},c}(Y)\to DM_{{\cyrrm{B}},c}(X)$. The functor $f^*$ is monoidal, i.e. $f^*\mathbbold{1}_X=\mathbbold{1}_Y$.
\item For any separated morphism of finite type $f:Y\to X$, there is a functor $f_!:DM_{{\cyrrm{B}},c}(Y)\to DM_{{\cyrrm{B}},c}(X)$, which has a right adjoint $f^!:DM_{{\cyrrm{B}},c}(X)\to DM_{{\cyrrm{B}},c}(Y)$.
\item There is a structure of a covariant (resp. contravariant) 2-functors on $f\mapsto f_*$ and $f\mapsto f_!$ (resp. $f\mapsto f^*$ and $f\mapsto f^!$).
\item For any proper morphism $f$, there is a natural identification $f_!=f_*$.
\item 
\label{rel_purity}
(Relative purity) For any smooth morphism $f$ of relative dimension $d$, there is an isomorphism of functors $P(f):f^!\simeq f^*(d)[2d]$.
\item
\label{BC}
 (Base change isomorphism) For any cartesian square of schemes
$$
  \xymatrix{
    Y' \ar[r]^-{f'} \ar[d]^-{g'} & X' \ar[d]^-{g}\\
    Y \ar[r]^-{f} & X,
  }
$$
such that $f$ is separated of finite type, there  is an isomorphism $Ex(\Delta^*_!):g^*f_!\simeq f'_!g'^*$.
\item 
\label{proj_formula}
(Projection formula) For any separated morphism of finite type $f:Y\to X$, there is an 
isomorphism $(f_!K)\otimes_XL\simeq f_!(K\otimes_Yf^*L)$.
\item 
\label{loc_dist_triangle}
(Localization distinguished triangle) For any closed immersion $i:Z\to S$ with complementary open immersion $j$, there is a distinguished triangle of natural transformations
$$
j_!j^!\to1\to i_*i^*\to j_!j^![1]
$$
where the first two maps are canonical adjunction maps.
\item 
\label{abs_purity}
(Absolute purity) For any closed immersion $i:Z\to S$ between regular schemes of codimension $c$, there is an isomorphism $P(i):i^!\mathbbold{1}_S\simeq \mathbbold{1}_Z(-c)[-2c]$.
\end{enumerate}
\end{recall}

\subsection{Borel-Moore motives and basic functorialities}\

\begin{definition}

Let $f:X\to S$ be a separated morphism of finite type. The \textbf{Borel-Moore motive} of $X$ relative to $S$, denoted by $M^{BM}(X/S)$, is defined as the element $f_!f^*\mathbbold{1}_S=f_!\mathbbold{1}_X$ in the category $DM_{{\cyrrm{B}},c}(S)$.

\end{definition}

\begin{remark}
The following properties follow from the definition above:
\begin{enumerate}
\item For any scheme $Y$, $M^{BM}(Y/Y)=\mathbbold{1}_Y$.

\item If we have a tower $X\to Y\overset{p}{\to}S$, then $M^{BM}(X/S)=p_!M^{BM}(X/Y)$. 
\end{enumerate}
The second point shows that for functorial compatibilities of Borel-Moore motives we may work over the base $Y$ and then "go down" over the base $S$.

\end{remark}

Borel-Moore motives have three basic functorialities:


\begin{lemma}
\label{bm_basic}
Let $p:Y\to S$ be a separated morphism of finite type.
\begin{enumerate}
\item
(Proper functoriality) 
Let $f:X\to Y$ be a proper morphism, then there is a map $f^*:M^{BM}(Y/S)\to M^{BM}(X/S)$.
\item
\label{sm_funct}
(Smooth functoriality) 
Let $f:X\to Y$ be a smooth morphism of relative dimension $d$, then there is a map $f_*:M^{BM}(X/S)\to M^{BM}(Y/S)(-d)[-2d]$.
\item 
\label{localization_leq}
(Localization triangle)
Let $Z$ be a closed subscheme of $Y$, then there is a canonical map $\partial_{Y,Z}:M^{BM}(Z/S)\to M^{BM}(Y-Z/S)[1]$, called boundary map, such that the triangle
$$
M^{BM}(Y-Z/S)\to M^{BM}(Y/S)\to M^{BM}(Y/S)\overset{\partial_{Y,Z}}{\to}M^{BM}(Y-Z/S)[1]
$$
is a distinguished triangle in $DM_{{\cyrrm{B}},c}(S)$, where the two first maps are the proper and smooth functorialities constructed previously. 

\end{enumerate}
\end{lemma}

\proof
\begin{enumerate}
\item The map is constructed as
$$
M^{BM}(Y/S)=p_!\mathbbold{1}_Y\overset{ad_{(f^*,f_*)}}{\longrightarrow} p_!f_*f^*\mathbbold{1}_Y=M^{BM}(X/S).
$$

\item Using the relative purity (Recall~\ref{recall_DM} \ref{rel_purity})) of $f$, we construct the map as
\begin{align*}
        &M^{BM}(X/S)=p_!f_!f^*\mathbbold{1}_Y\overset{(P(f))^{-1}}{\simeq} p_!f_!f^!\mathbbold{1}_Y(-d)[-2d]\\
        &\overset{ad'_{(f_!,f^!)}}{\longrightarrow} p_!\mathbbold{1}_Y(-d)[-2d]=M^{BM}(Y/S)(-d)[-2d].
\end{align*}

\item The localization triangle comes from the localization distinguished triangle (Recall~\ref{recall_DM} (\ref{loc_dist_triangle})).
\end{enumerate}
\endproof

\begin{remark}
\label{BM_union}
\begin{enumerate}
\item The proper and smooth functorialities are compatible with compositions. The proper case is clear, and the smooth case is stated in Recall~\ref{purity_comp} \ref{composition_rel_purity}).

\item 
\label{BM_union_1}
By the localization triangle, we know that if $X=X_1\coprod X_2$ is a disjoint union of two sub-$S$-schemes, then there is a canonical identification $M^{BM}(X/S)= M^{BM}(X_1/S)\oplus M^{BM}(X_2/S)$. Via this identification, the proper functoriality and the boundary map are additive (\cite[Proposition 2.26]{Deg2}): for example in the case of the proper functoriality, if $X=X_1\coprod X_2$ and $Y=Y_1\coprod Y_2$ are disjoint unions, and $f:X\to Y$ is a proper morphism obtained by gluing two proper morphisms $f_1:X_1\to Y_1$ and $f_2:X_2\to Y_2$, then $f_*=f_{1*}+f_{2*}$.

\end{enumerate}
\end{remark}

In what follows, we are going to check that some other operations are compatible with the three basic functorialities. The following lemma is a general fact in triangulated categories:
\begin{lemma}
\label{uniqueness_h}
Suppose we have the following diagram in a triangulated category
$$
  \xymatrix{
    A \ar[r]^-{u} \ar[d]^-{f} \ar@{}[rd]|{(1)} & B \ar[r]^-{v} \ar[d]^-{g} & C \ar[r]^-{w} \ar@{.>}[d]^-{h} & A[1]\ar[d]^-{f[1]}\\
    A'  \ar[r]^-{u'} & B \ar[r]^-{v'} & C' \ar[r]^-{w'} & A'[1],
  }
$$
where the two rows $(u,v,w)$ and $(u',v',w')$ are distinguished triangles and the maps $f$ and $g$ are such that square $(1)$ commutes. Suppose in addition that $Hom(A[1],C')=0$. Then there exists one and only one map $h:C\to C'$ which makes the whole diagram commutative.
\end{lemma}
In particular, by \cite[Proposition 2.3.3]{CD}, the conditions in Lemma~\ref{uniqueness_h} are satisfied for the localization triangle, which means that in order to check the compatibility between boundary maps, it is sufficient to check the compatibility between proper and smooth functoriality maps.

Now we check the basic functorialities are compatible one with each other. Until the end of this section, all schemes are supposed to be separated of finite type over a base scheme $S$. The key point is the following lemma, namely the compatibility between proper and smooth functorialities for a cartesian square:
\begin{lemma}
\label{BC_proper_smooth}
Suppose we have a cartesian square
$$
  \xymatrix{
    Y' \ar[r]^-{q} \ar[d]^-{g} & X' \ar[d]^-{f} \\
    Y \ar[r]^-{p} & X,
  }
$$
with $p$ a proper morphism and $f$ a smooth morphism of relative dimension $d$. Then we have a commutative diagram
$$
  \xymatrix{
    M^{BM}(Y'/S)\ar[d]^-{g_*} & M^{BM}(X'/S)  \ar[l]_-{q^*} \ar[d]^-{f_*} \\
    M^{BM}(Y/S)(-d)[-2d]  & M^{BM}(X/S)(-d)[-2d]. \ar[l]_-{p^*} 
  }
$$
\end{lemma}
The proof of Lemma~\ref{BC_proper_smooth} will be given in Section~\ref{proof_BC_proper_smooth}. The three following lemmas are straightforward consequences of Lemma~\ref{uniqueness_h} and Lemma~\ref{BC_proper_smooth}:
\begin{lemma}
\label{cartesian_diff_comp}
Suppose we have a cartesian square
$$
  \xymatrix{
    Z' \ar[r]^-{i'} \ar[d]^-{q} & X' \ar[d]^-{p} \\
    Z \ar[r]^-{i} & X,
  }
$$
with $i$ a closed immersion and $p$ a separated morphism of finite type, which can be completed to a diagram of two cartesian squares
$$
  \xymatrix{
    Z' \ar[r]^-{i'} \ar[d]^-{q} & X' \ar[d]^-{p}& X'-Z' \ar[l]_-{j'} \ar[d]^-{r} \\
    Z \ar[r]^-{i} & X & X-Z. \ar[l]_-{j}
  }
$$
\begin{enumerate}
\item If $p$ is a proper morphism, then we have a commutative diagram
$$
  \xymatrix{
    M^{BM}(X'-Z'/S)[1] & M^{BM}(Z'/S)  \ar[l]_-{\ \ \ \ \partial_{X',Z'}} \\
    M^{BM}(X-Z/S)[1] \ar[u]^-{r^*} & M^{BM}(Z/S). \ar[l]_-{\ \ \ \partial_{X,Z}} \ar[u]^-{q^*}
  }
$$

\item If $p$ is a smooth morphism of relative dimension $d$, then we have a commutative diagram
$$
  \xymatrix{
    M^{BM}(X'-Z'/S)(d)[2d+1] \ar[d]^-{r_*} & M^{BM}(Z'/S)(d)[2d]  \ar[l]_-{\ \ \ \ \partial_{X',Z'}} \ar[d]^-{q_*} \\
    M^{BM}(X-Z/S)[1]  & M^{BM}(Z/S). \ar[l]_-{\ \ \ \partial_{X,Z}} 
  }
$$
\end{enumerate}
\end{lemma}

\begin{lemma}
\label{comp_diff}
Suppose we have two consecutive closed immersions $W\to Y\to Z$. Then we have a commutative diagram
$$
  \xymatrix{
    M^{BM}(Z-Y/S) [1]  & & M^{BM}(Y-W/S). \ar[ll]_-{\partial_{Z-W,Y-W}} \ar[ld]^-{j_*} \\
     & M^{BM}(Y/S) \ar[lu]^-{\partial_{Z,Y}} &
  }
$$
\end{lemma}

\proof
The result follows from Lemma~\ref{cartesian_diff_comp} $(2)$ applied to the 
following commutative diagram 
$$
  \xymatrix{
    Y-W \ar[r]^-{i_2}\ar@{}[rd]|{(1)} \ar[d]^-{j_3} & Z-W \ar[d]^-{j_4} &\\
    Y \ar[r]_-{i_1} & Z & Z-Y, \ar[l]^-{j_1} \ar[lu]_-{j_2}
  }
$$
where the $i_k$ are closed immersions and the $j_k$ are open immersions, and the square $(1)$ is cartesian. 

\endproof

\begin{lemma}
\label{compat_diff_two}

Suppose that we have the following commutative diagram of schemes:
$$
  \xymatrix{
    Z' \ar[r]^-{i_2} \ar[d]_-{i_3}& X' \ar[d]_-{i_4} &\\
    Z \ar[r]^-{i_1} & X & Y, \ar[l]_-{j} \ar[lu]_-{j'}
  }
$$
where $j$ and $j'$ are open immersions and all $i_k$'s are closed immersions, such that $X-Z=X'-Z'=Y$. Then the following diagram is commutative:
$$
  \xymatrix{
    M^{BM}(Y/S) [1]  & & M^{BM}(Z/S). \ar[ll]_-{\partial_{X,Z}} \ar[ld]^-{i_3^*} \\
     & M^{BM}(Z'/S) \ar[lu]^-{\partial_{X',Z'}} &
  }
$$
\end{lemma}

\subsection{Purity isomorphism, Gysin morphism and refined Gysin morphism}\
\label{chapter_gysin}

\noindent
In this section, $S$ is a base scheme, and all $S$-schemes 
are supposed to be separated of finite type.

\begin{remark}
\label{remark_qproj}
(\cite[Remarks 1.1.3]{Lev1})) If $f:X\to Y$ is a quasi-projective morphism, $f$ is always \textbf{smoothly embeddable}, i.e. $f$ can be factorized as $X\overset{i}{\to}P\overset{p}{\to}Y$, where $i$ is closed immersion and $p$ is a smooth morphism.
Furthermore, given two quasi-projective morphisms $f:X\to Y$ and $g:Y\to Z$, their composition is quasi-projective, thus smoothly embeddable. In addition, there exists a commutative diagram
$$
  \xymatrix{
     X \ar[r]^-{i_1} \ar[rd]_-{f} & P \ar[r]^-{i_3} \ar[d]_-{p_1} \ar@{}[rd]|{(1)} & R \ar[d]^-{p_3}\\
       &  Y  \ar[r]^-{i_2}  \ar[rd]_-{g} & Q \ar[d]^-{p_2} \\
        &  & Z,
  }
$$
where all $i_k$'s are closed immersions and $p_k$'s are smooth morphisms, and the square $(1)$ is cartesian.
\end{remark}

\begin{definition}
\label{relative_dim}
Let $f:X\to Y$ be a quasi-projective morphism between regular schemes, which is factorized as the composite map $X\overset{i}{\to}P\overset{p}{\to}Y$, where $i$ is a closed immersion and $p$ is a smooth morphism. The scheme $P$ is then regular, and therefore $i$ is a regular closed immersion (\cite[B.7.1]{Ful}). If the immersion $X\overset{i}{\to}P$ is of codimension $d_1$ (i.e. its normal bundle is of rank $d_1$ over $X$) and the smooth morphism $P\overset{p}{\to}Y$ is of relative dimension $d_2$, we say that the morphism $f$ is of relative dimension $d=d_2-d_1$.
\end{definition}

Note that the relative dimension is always Zariski locally well defined and independent of the factorization. Recall that there are two types of purity isomorphisms: the relative purity (Recall~\ref{recall_DM} \ref{rel_purity})) and the absolute purity (Recall~\ref{recall_DM} \ref{abs_purity})).
Both of them are stable by composition:
\begin{recall}
\label{purity_comp}
\begin{enumerate}
\item
\label{composition_rel_purity}
(\cite[Remark 2.4.52]{CD}) Let $X\overset{f}{\to}Y\overset{g}{\to}Z$ be two consecutive smooth morphisms. Then $P(g\circ f)=P(g)\circ P(f)$.
\item (\cite[Theorem 2.4.9]{Deg5}) Let $X\overset{i_1}{\to}Y\overset{i_2}{\to}Z$ be two consecutive closed immersions between regular schemes. Then $P(i_2\circ i_1)=P(i_2)\circ P(i_1)$.
\end{enumerate}
\end{recall}

Furthermore, the two kind of purity isomorphisms are compatible with each other in the following way:
\begin{lemma}
\label{compat_purity}
Consider a cartesian square of regular schemes
$$
  \xymatrix{
    Z' \ar[r]^-{i'} \ar[d]_-{g} & X' \ar[d]^f \\
    Z \ar[r]^-{i} & X, 
  }
$$
where $i$ and $i'$ are closed immersions of codimension $c$, $f$ and $g$ are smooth morphisms of relative dimension $d$. Then we have $P(i')\circ P(f)=P(g)\circ P(i)$.
\end{lemma}
The proof of Lemma~\ref{compat_purity} will be given in Section~\ref{proof_BC_proper_smooth}. Combining the relative purity and absolute purity, we generalize the purity isomorphism to any quasi-projective morphism between regular schemes:
\begin{definition}
\label{def_purity}
Let $f:X\to Y$ be a quasi-projective morphism between regular schemes of relative dimension $d$. 
We have the following \textbf{purity isomorphism}:
$$
P(f):f^!\mathbbold{1}_Y\simeq f^*\mathbbold{1}_Y(d)[2d]=\mathbbold{1}_X(d)[2d].
$$

\end{definition}

Firstly we need to check that 
the purity isomorphism is independent of the choice of factorization. We need the following lemma:
\begin{lemma}
\label{compose_purity}
Consider two consecutive morphisms between regular schemes $X\overset{i}{\to}P\overset{p}{\to}Y$, where $i$ is a closed immersion, $p$ is a smooth morphism, such that $p\circ i$ is a closed immersion. Then we have $P(p\circ i)=P(p)\circ P(i)$.

\end{lemma}
\proof
We have the following commutative diagram of schemes:
$$
  \xymatrix{
   & P \ar[r]^-{p} \ar@{}[rd]|{(1)} & Y \\
   & P\times_SX \ar[r]^-{p_2} \ar[u]_-{p_1} & X, \ar[u]_-{p\circ i}  \\
    X \ar[ruu]^-{i} \ar[ru]^-{\alpha} \ar[rru]_-{id_X} & &
  }
$$
where $\alpha=(i,id_X):X\to P\times_SX$ is a closed immersion, as well as $p_1$. Then Recall~\ref{purity_comp} implies that $P(i)=P(p_1)\circ P(\alpha)$. Also, applying Lemma~\ref{compat_purity} to the cartesian square $(1)$, we have $P(p)\circ P(p_1)=P(p\circ i)\circ P(p_2)$. Therefore in order to show $P(p\circ i)=P(p)\circ P(i)$, we only need to show that $P(p_2)\circ P(\alpha)=P(id_X)$, i.e. we are reduced to the case where $X=Y$ and $i$ is a section of $p$, which is proved in \cite{Deg5}.

\endproof

\begin{corollary}
\label{purity_indep}
\begin{enumerate}
\item The purity isomorphism does not depend on the factorization.
\item 
\label{purity_composition}
The purity isomorphism is compatible with composition.
\end{enumerate}
\end{corollary}

\proof

\begin{enumerate}
\item (\cite[Proposition 6.6(a)]{Ful}) Let $f:X\to Y$ be a quasi-projective morphism and suppose that there are two factorizations $X\to P\to Y$ and $X\to P'\to Y$ of $f$. Compare both of them to the diagonal
$$
  \xymatrix{
  X \ar[r] & P\times_YP' \ar[r] \ar[d] & P' \ar[d] \\
    & P \ar[r] & Y,
  }
$$
and the corollary follows from Lemma~\ref{compose_purity}.
\item Results from Remark~\ref{remark_qproj}, Recall~\ref{purity_comp} and Lemma~\ref{compat_purity}.
\end{enumerate}
\endproof

\begin{definition}
\label{def_gysin}
Consider two consecutive morphisms $X\overset{f}{\to}Y\overset{p}{\to}S$, where $X$ and $Y$ are regular $S$-schemes and $f$ is a quasi-projective morphism of relative dimension $d$. 
We deduce from the purity isomorphism a \textbf{Gysin morphism}
\begin{align*}
        G(f):&M^{BM}(X/S)=p_!f_!\mathbbold{1}_X\overset{P(f)^{-1}}{\simeq} p_!f_!f^!\mathbbold{1}_Y(-d)[-2d]\\
        &\xrightarrow{ad'_{(f_!,f^!)}} p_!\mathbbold{1}_Y(-d)[-2d]=M^{BM}(Y/S)(-d)[-2d].
\end{align*}
\end{definition}

It follows from Corollary~\ref{purity_indep} \ref{purity_composition}) that the formation of Gysin morphisms is compatible with composition:
\begin{corollary}
\label{comp_gysin}
If $X\overset{f}{\to}Y\overset{g}{\to}Z$ are two quasi-projective morphisms between regular schemes, then $G(g\circ f)=G(g)\circ G(f)$.

\end{corollary}

Now we check the compatibility between Gysin morphisms with the three basic functorialities. In particular, if $f$ is a smooth morphism between regular schemes, then $G(f)$ is nothing but the smooth functoriality (Lemma~\ref{bm_basic} \ref{sm_funct})) for Borel-Moore motives, so by Corollary~\ref{comp_gysin}, 
all Gysin morphisms are compatible with the smooth functoriality. For the proper functoriality, it is compatible with Gysin morphisms in the transversal intersection case:

\begin{lemma}
\label{compat_gysin_func}
Consider a cartesian square of regular $S$-schemes
$$
  \xymatrix{
    Z' \ar[d]_-{q} \ar[r]^-{i'} & X' \ar[d]^-{p} \\ 
    Z \ar[r]^i & X,
  }
$$
where $p$ and $q$ are proper morphisms and $i$ and $i'$ are closed immersions, such that $p$ is transversal to $i$. 
\footnote{Recall that this means that the canonical closed immersion from the normal cone of $Z'$ in $X'$ to the pull-back by $p$ of the normal cone of $Z$ in $X$ (\cite[B.6.1]{Ful}) is an isomorphism.}
Then the following square is commutative:
$$
  \xymatrix{
    M^{BM}(Z/S) \ar[d]_-{q^*} \ar[r]^-{G(i)}& M^{BM}(X/S)(c)[2c]  \ar[d]^-{p^*} \\
    M^{BM}(Z'/S) \ar[r]^-{G(i')} & M^{BM}(X'/S)(c)[2c].
  }
$$

In the case where $p$ and $q$ are closed immersions, then the following square is commutative:
$$
  \xymatrix{
    M^{BM}(Z/S) \ar[d]_-{\partial_{Z,Z'}} \ar[r]^-{G(i)}& M^{BM}(X/S)(c)[2c]  \ar[d]^-{\partial_{X,X'}} \\
    M^{BM}(Z-Z'/S)[1] \ar[r]^-{G(i'')} & M^{BM}(X-X'/S)(c)[2c+1].
  }
$$
\end{lemma}
The proof of Lemma~\ref{compat_gysin_func} will be given in Section~\ref{proof_BC_proper_smooth}.
\medskip

Now we introduce a functorial base change property of Borel-Moore motives, which comes from the base change isomorphism (Recall~\ref{recall_DM} \ref{BC})):
\begin{definition}
\label{def_monoidal_BC}
Consider a cartesian square of schemes
$$
  \xymatrix{
    X' \ar[r]^-{q} \ar[d]_-{g} \ar@{}[rd]|{\Delta} & Y' \ar[d]^-{f} \\
    X \ar[r]^-{p} & Y,
  }
$$
where $p$ and $q$ are separated of finite type. Then there is an isomorphism
$$
BC(\Delta):M^{BM}(X'/Y')
=
q_!g^*\mathbbold{1}_{X}
\overset{(Ex(\Delta^*_!))^{-1}}{\simeq} 
f^*p_!\mathbbold{1}_{X}
=
f^*M^{BM}(X/Y)
$$
called \textbf{base change isomorphism} for Borel-Moore motives.

\end{definition}

\begin{remark}
\label{BC_transitive}
It follows from \cite[Corollary 2.2.12]{CD} that the base change isomorphism is compatible with horizontal and vertical compositions of squares.
\end{remark}

As an application, we deduce from the base change isomorphism a formula of K\"unneth type for Borel-Moore motives:
\begin{lemma}[K\"unneth formula]
\label{kunneth}
For any $S$-schemes $X$ and $Y$, there is an isomorphism $M^{BM}(X\times_SY/S)\simeq M^{BM}(X/S)\otimes M^{BM}(Y/S)$.

\end{lemma}

\proof

Denote by $f:X\to S$ the structural morphism. Then using the projection formula (Recall~\ref{recall_DM} \ref{proj_formula})) and the base change isomorphism, we have
\begin{align*}
 M^{BM}(X/S)\otimes M^{BM}(Y/S) &= f_!\mathbbold{1}_X\otimes M^{BM}(Y/S) \simeq f_!(\mathbbold{1}_X\otimes f^*M^{BM}(Y/S))\\
 &\simeq f_!(\mathbbold{1}_X\otimes M^{BM}(X\times_SY/X))=M^{BM}(X\times_SY/S).
\end{align*}
\endproof

\begin{definition}
\label{refine_cartesian}
Take the notations in Definition~\ref{def_monoidal_BC}. 
Let $\phi:M^{BM}(X/Y)\to \mathbbold{1}_{Y}(m)[n]$ be a map, and let $r:Y'\to S$ be a separated morphism of finite type. We define a map $R_f(\phi/S):M^{BM}(X'/S)\to M^{BM}(Y'/S)(m)[n]$ as the composition 
\begin{align*}
        &M^{BM}(X'/S)=r_!M^{BM}(X'/Y')\overset{BC(\Delta)}{\simeq}r_!f^*M^{BM}(X/Y)\\
        &\overset{\phi}{\longrightarrow}r_!f^*\mathbbold{1}_{Y}(m)[n]=M^{BM}(Y'/S)(m)[n].
\end{align*}
Similarly, if $\phi:\mathbbold{1}_{Y}\to M^{BM}(X/Y)(m)[n]$ is a map, we define a map $$
R_f(\phi/S):M^{BM}(Y'/S)\to M^{BM}(X'/S)(m)[n].
$$
\end{definition}
It follows from Remark~\ref{BC_transitive} that the $R$-operation is transitive:
\begin{lemma}
\label{compat_r}
For any $Y''\overset{g}{\to}Y'\overset{f}{\to}Y$ we have $R_{f\circ g/S}(\phi)=R_g(R_f(\phi/S)/S)$.

\end{lemma}
When $f$ and $S$ are clear, we take also the notation $R_f(\phi)$ or $R(\phi)$ instead of $R_f(\phi/S)$ for simplicity.
Now we check that the $R$-operation is compatible with basic functorialities:
\begin{lemma}
\label{lemma_monoidal_BC}
Take the notations in Definition~\ref{def_monoidal_BC}. Then
\begin{enumerate}
\item
\label{lemma_monoidal_BC1} 
If $p$ and $q$ are proper morphisms, then $R_f(p^*)=q^*$. In other words, the following diagram is commutative:
$$
  \xymatrix{
    \mathbbold{1}_{Y'} \ar@{=}[d] \ar[r]^-{q^*} & M^{BM}(X'/Y') \ar[d]^-{BC(\Delta)}_-{\wr} \\
    f^*\mathbbold{1}_{Y} \ar[r]^-{p^*} &  f^*M^{BM}(X/Y).
  }
$$

\item If $p$ and $q$ are smooth morphisms of relative dimension $d$, then $R_f(p_*)=q_*$. In other words, the following diagram is commutative:
$$
  \xymatrix{
    M^{BM}(X'/Y') \ar[d]^-{BC(\Delta)}_-{\wr} \ar[r]^-{q_*} &  \mathbbold{1}_{Y'}(-d)[-2d] \ar@{=}[d] \\
    f^*M^{BM}(X/Y) \ar[r]^-{p_*} &    f^*\mathbbold{1}_{Y}(-d)[-2d].
  }
$$

\item If $p$ and $q$ are closed immersions, then the following diagram is commutative:
$$
  \xymatrix{
    M^{BM}(X'/Y') \ar[d]^-{BC(\Delta)}_-{\wr} \ar[r]^-{\partial_{Y',X'}} & M^{BM}(Y'-X'/Y')[1] \ar[d]^-{BC(\Delta')}_-{\wr} \\
    f^*M^{BM}(X/Y) \ar[r]^-{\partial_{Y,X}} & f^*M^{BM}(Y-X/Y)[1],
  }
$$
where $\Delta'$ is the cartesian square
$$
  \xymatrix{
    Y'-X' \ar[r]^-{q_1} \ar[d]^-{g} \ar@{}[rd]|-{\Delta'} & Y' \ar[d]^-{f} \\
    Y-X \ar[r]^-{p_1} & Y.
  }
$$
\end{enumerate}
\end{lemma}
The proof of Lemma~\ref{compat_gysin_func} will be given in Section~\ref{proof_BC_proper_smooth}.
\medskip

For any scheme $X$, the boundary map $\mathbbold{1}_X\overset{\partial_{\mathbb{A}^1_X,X}}{\longrightarrow}M^{BM}(\mathbb{G}_{m,X}/X)[1]$ associated to the zero section $X\to\mathbb{A}^1_X$ has a canonical right inverse
$$
\phi_X:M^{BM}(\mathbb{G}_{m,X}/X)[1]\to\mathbbold{1}_X,
$$
such that for any scheme $Y$ over $X$, $R(\phi_X)=\phi_Y$.
\footnote{$\phi_X$ is an element of $H^{BM}_{1,0}(\mathbb{G}_{m,X}/X)
\simeq \mathbb{Z}\times \mathcal{O}^*(X)$, which corresponds to the fundamental parameter of $\mathbb{G}_m$, see Lemma~\ref{elem_BM}.}
It results from the construction of the K\"unneth formula (Lemma~\ref{kunneth}) that the map $\phi_X$ is compatible with tensor products:

\begin{lemma}
\label{compat_phi_tensor}
Let $f:X\to S$ be a separated morphism of finite type. Then the following diagram is commutative:
$$
  \xymatrix{
    M^{BM}(\mathbb{G}_{m,X}/S) \ar[d]_-{\wr} \ar[r]^-{f_!\phi_X} & M^{BM}(X/S).\\
    M^{BM}(X/S)\otimes M^{BM}(\mathbb{G}_{m,S}/S)[1] \ar[ur]_-{id\otimes\phi_S} & 
  }
$$

\end{lemma}

\begin{definition}
Let $X$ be a regular $S$-scheme endowed with a $S$-morphism $f:X\to\mathbb{G}_{m,S}$. Then the morphism $\Gamma_f:X\to X\times\mathbb{G}_{m,S}=\mathbb{G}_{m,X}$, the graph of $f$, is a closed immersion of codimension $1$ between regular schemes. Define the map $\psi(f):M^{BM}(X/S)[1]\to M^{BM}(X/S)(1)[2]$ as the composition 
$$
M^{BM}(X/S)[1]\overset{G(\Gamma_f)}{\longrightarrow} M^{BM}(\mathbb{G}_{m,X}/S)(1)[3]\overset{\phi_X}{\longrightarrow} M^{BM}(X/S)(1)[2].
$$

\end{definition}

\begin{proposition}
\label{double_deformation}
. 
Let $X$ be a regular $S$-scheme endowed with a flat morphism $X\to\mathbb{A}^1$, such that the fiber over the point $0$ is regular, denoted by $Z$. Denote by $f:X-Z\to\mathbb{G}_m$, $i:Z\to X$ and $j:X-Z\to X$ the corresponding morphisms. Then the following diagram is commutative:
$$
  \xymatrix{
M^{BM}(Z/S) \ar[r]^-{\partial_{X,Z}} \ar[rd]_-{G(i)}& M^{BM}(X-Z/S)[1] \ar[r]^-{\psi(f)} \ar@{}[d]|-{(1)} & M^{BM}(X-Z/S)(1)[2]. \ar[ld]^-{j_*} \\
 & M^{BM}(X/S)(1)[2] &
  }
$$
\end{proposition}

\proof

We use the same strategy as in \cite[Proposition 2.6.5]{Deg4}. Without loss of generality, suppose that the closed immersion $i:Z\to X$ is of codimension $c$.

Compare the diagram $(1)$ to the following diagram
$$
  \xymatrix{
M^{BM}(N_ZX/S)(c)[2c] \ar[r]^-{\partial_{D_ZX,N_ZX}} \ar[rd]_-{G(i_1)}& M^{BM}(\mathbb{G}_{m,X}/S)(c)[2c+1] \ar[r]^-{\psi(f_1)} \ar@{}[d]|-{(2)} & M^{BM}(\mathbb{G}_{m,X}/S)(c+1)[2c+2]. \ar[ld]^-{j_{1*}} \\
 & M^{BM}(D_ZX/S)(c+1)[2c+2] &
  }
$$
There is a natural map from the diagram $(1)$ to the diagram $(2)$ given by termwise Gysin morphisms. By Corollary~\ref{comp_gysin} and Lemma~\ref{compat_gysin_func}, all the maps in the diagram are compatible with Gysin morphisms (the flatness hypothesis shows that the intersection is transversal). In addition, the closed immersion $X\to D_ZX$ is a section of the structure morphism $D_ZX\to X$.  Consequently, to show that the diagram $(1)$ commutes, it is sufficient to show that the diagram $(2)$ commutes.

Now using the same strategy, compare the diagram $(2)$ to the following diagram
$$
  \xymatrix{
M^{BM}(Z/S) \ar[r]^-{\partial_{\mathbb{A}^1_Z,Z}} \ar[rd]_-{G(i_2)}& M^{BM}(\mathbb{G}_{m,Z}/S)[1] \ar[r]^-{\psi(f_2)} \ar@{}[d]|-{(3)} & M^{BM}(\mathbb{G}_{m,Z}/S)(1)[2], \ar[ld]^-{j_{2*}} \\
 & M^{BM}(\mathbb{A}^1_Z/S)(1)[2] &
  }
$$
and we know that it suffices to show the commutativity of the diagram $(3)$, which follows from Lemma~\ref{compat_phi_tensor} and the fact that the map $\partial_{\mathbb{A}^1_Z,Z}$ is a section of $\phi_Z$.

\endproof

Now we are going to apply Proposition~\ref{double_deformation} to the framework of the deformation to the normal cone (\cite[5.1]{Ful}). 
Consider a closed immersion $i:X\to Y$ of codimension $c$ between regular schemes. Let $N_XY$ be the normal cone of $X$ in $Y$, and $D_XY$ the deformation space, which is the difference $Bl_{X\times \{0\}}(Y\times\mathbb{A}^1)-Bl_XY$ (\emph{ibid.}, where $\mathbb{P}^1$ is replaced by $\mathbb{A}^1$ here). The scheme $D_XY$ is a scheme over $Y\times\mathbb{A}^1$, which is flat over $\mathbb{A}^1$. The fiber over the point $0$ is the normal cone $N_XY$, and the complementary $D_XY-N_XY$ is isomorphic to $\mathbb{G}_{m,Y}$. We have the following commutative diagram:
$$
  \xymatrix{
    N_XY \ar[r]^-{i_2} \ar[d]^-{p} & D_XY \ar[d]^-{f} & \mathbb{G}_{m,Y}. \ar[l]_-{j} \ar[ld]^-{\pi_1} \\
    X \ar[r]^-{i_1} & Y &
  }
$$
The maps $i_1$ and $i_2$ are closed immersions and $j$ is the complementary open immersion to $i_2$. Denote by $\pi_2$ the projection to the second factor $Y\times \mathbb{G}_m\to \mathbb{G}_m$.
\begin{proposition}
\label{compat_gysin_diff}
The following diagram is commutative:
$$
  \xymatrix{
    M^{BM}(N_XY/Y) \ar[r]^-{\partial_{D_XY,N_XY}} \ar[d]^-{p_*=G(p)}& M^{BM}(\mathbb{G}_{m,Y}/Y)[1] \ar[d]^-{\phi_Y} \\
    M^{BM}(X/Y)(-c)[-2c] \ar[r]^-{G(i_1)} & \mathbbold{1}_Y.
  }
$$
\end{proposition}

\proof
Apply Proposition~\ref{double_deformation} to the quasi-projective $Y$-scheme $D_XY$. By Corollary~\ref{comp_gysin}, we obtain the following commutative diagram:  
$$
  \xymatrix{
    M^{BM}(N_XY/Y) \ar[r]^-{\partial_{D_XY,N_XY}} \ar[rd]^-{G(i_2)} \ar[rdd]_-{G(i_1\circ p)} & M^{BM}(\mathbb{G}_{m,Y}/Y)[1] \ar[r]^-{\psi(\pi_2)} & M^{BM}(\mathbb{G}_{m,Y}/Y)(1)[2]. \ar[ld]_-{j_*} \ar[ldd]^-{\pi_{1*}}\\
    & M^{BM}(D_XY/Y)(1)[2] \ar[d]^-{G(f)} & \\
    & \mathbbold{1}_Y &
  }
$$
Therefore the proposition follows from the following commutative diagram:
$$
  \xymatrix{
    M^{BM}(\mathbb{G}_{m,Y}/Y)[1] \ar[r]^-{\psi(\pi_2)} \ar[d]_-{\phi_Y} & M^{BM}(\mathbb{G}_{m,Y}/Y)(1)[2], \ar[ld]^-{\pi_{1*}} \\
    \mathbbold{1}_Y &
  }
$$
which follows from Lemma~\ref{compat_phi_tensor}.
\endproof

\begin{definition}
\label{refined_gysin}
Let $i:X\to Y$ be a closed immersion of codimension $c$ between regular schemes. 
Let $f:Y'\to Y$ be a morphism such that $Y'$ is a $S$-scheme, and form the cartesian square
$$
  \xymatrix{
    X' \ar[r]^-{i'} \ar[d]_-{g}  \ar@{}[rd]|-{\Delta} & Y' \ar[d]^-{f} \\ 
    X \ar[r]^-{i} & Y.
  }
$$
The map $R_f(i)=R_f(G(i)/S):M^{BM}(X'/S)\to M^{BM}(Y'/S)(c)[2c]$ is called the \textbf{refined Gysin morphism} associated to the square $\Delta$.
\end{definition}


\begin{corollary}
\label{rgysin_1}
Take the setting of Definition~\ref{refined_gysin}. Let $N$ be the normal cone of $X$ in $Y$, and let $N'$ be the fiber product $N\times_YY'$. Then the following diagram is commutative:
$$
  \xymatrix{
    M^{BM}(N'/S) \ar[r]^-{\partial_{D_ZX,N_ZX}} \ar[d]_-{p_*}& M^{BM}(\mathbb{G}_{m,Y'}/S)[1] \ar[d]^-{\phi_{Y'}} \\
    M^{BM}(X'/S)(-c)[-2c] \ar[r]^-{R_f(i)} & M^{BM}(Y'/S).
  }
$$
\end{corollary}

\proof

Apply the R-operation to the diagram in Proposition~\ref{compat_gysin_diff}, then use Lemma~\ref{lemma_monoidal_BC} to show that the corresponding maps are the same.

\endproof

\begin{definition}
\label{def_specialization}
Let 
$i:V\to W$ be a closed immersion, and let $C_VW$ be the normal cone of $V$ in $W$, $D_VW$ be the deformation space. 
We define the following map, called \textbf{specialization map}.
$$
sp(i):M^{BM}(C_VW/S)\xrightarrow{\partial_{D_VW,C_VW}}M^{BM}(\mathbb{G}_{m,W}/S)[1]\xrightarrow{\phi_W}M^{BM}(W/S).
$$

\end{definition}

\begin{corollary}
\label{refined_gysin_compat}
Take the setting of Definition~\ref{refined_gysin}. Let $C'$ be the normal cone of $X'$ in $Y'$. Then the following diagram is commutative:
$$
  \xymatrix{
    M^{BM}(N'/S) \ar[r]^-{i^*} \ar[d]^-{p_*}& M^{BM}(C'/S) \ar[d]^-{sp(C'/Y')} \\
    M^{BM}(X'/S)(-c)[-2c] \ar[r]^-{R_f(i)} & M^{BM}(Y'/S).
  }
$$
\end{corollary}

\proof

The result follows from Corollary~\ref{rgysin_1} and Lemma~\ref{compat_diff_two}.


\endproof

\subsection{Composition laws}\
\medskip

\noindent
\textbf{Notations.}
In this section, $B$ is a regular scheme and $S$ a quasi-projective $B$-scheme. 
\begin{itemize}
\item For any $S$-schemes $X$ and $Y$, we denote by $XY=X\times_SY$ the fiber product over $S$, and $X\times Y=X\times_BY$ the fiber product over the base $B$. 
\item If $f:X\to B$ is a quasi-projective morphism with $X$ a regular scheme, denote by $d_X$ the relative dimension of $X$ over $B$ (Definition~\ref{relative_dim}) and denote by $G(X/B)=G(f):M^{BM}(X/B)\to \mathbbold{1}_B(-d_X)[-2d_X]$ the Gysin morphism associated to $f$ (Definition~\ref{def_gysin}).
\item For any morphism $X\to Y$, denote by $\delta_{X/Y}:X\to X\times_YX$ its diagonal embedding.
\end{itemize}

\begin{lemma}
\label{compat_BC_diag}
Let $p:X\to S$ be a proper morphism, with $\delta:X\to XX$ the diagonal embedding. Then the following diagram is commutative:
$$
  \xymatrix{
    M^{BM}(XX/S) \ar[r]^-{\delta^*} \ar[rd]^-{\sim}_{BC} & M^{BM}(X/S) \\
   & p_*p^*M^{BM}(X/S). \ar[u]_-{ad'_{(p^*,p_*)}}
  }
$$
\end{lemma}

\proof
The proof is a simple diagram chase of unit and counit maps, using the fact that diagonal embedding $\delta$ is a section of the projection to the second factor $XX\to X$.
\endproof

\begin{definition}
\label{BM_homology}
Let $X$ be a separated $S$-scheme of finite type. The \textbf{Borel-Moore homology} of $X$ relative to $S$ is defined as
$$
H^{BM}_{m,n}(X/S):=Hom_{DM_{{\cyrrm{B}},c}(S)}(M^{BM}(X/S),\mathbbold{1}_S(-n)[-m]).
$$
\end{definition}

\begin{lemma}
\label{Hom_BM}
Let $X$ and $Y$ be proper $S$-schemes that are smooth over $B$. 
Then there is an isomorphism
$$
\epsilon_{X,Y}:Hom_{DM_{{\cyrrm{B}},c}(S)}(M^{BM}(X/S),M^{BM}(Y/S)(q)[p])\simeq H^{BM}_{2d_Y-p,d_Y-q}(XY/B).
$$
\end{lemma}
\proof

Denote by $g:Y\to S$ and $h:Y\to B$ the structural morphisms. The isomorphism is constructed in the following way:
\begin{align*}
Hom(M^{BM}(X/S),M^{BM}(Y/S)(q)[p])
&=Hom(M^{BM}(X/S),g_*\mathbbold{1}_Y(q)[p]) \\
&=Hom(g^*M^{BM}(X/S),\mathbbold{1}_Y(q)[p]) \\
&\overset{BC}{\simeq} Hom(M^{BM}(XY/Y),\mathbbold{1}_Y(q)[p]) \\
&\overset{P(h)}{\simeq}Hom(M^{BM}(XY/Y),h^!\mathbbold{1}_B(q-d_Y)[p-2d_Y]) \\
&=Hom(M^{BM}(XY/B),\mathbbold{1}_B(q-d_Y)[p-2d_Y])\\
&=H^{BM}_{2d_Y-p,d_Y-q}(XY/B),
\end{align*}
where the first isomorphism is the inverse of the base change isomorphism for Borel-Moore motives (Definition~\ref{def_monoidal_BC}), and $P(h)$ is the purity isomorphism for the smooth morphism $h$.
\endproof

\begin{lemma}
\label{lemma_epsilon}
Take the setting in Lemma~\ref{Hom_BM}. Let $\alpha$ be a map from $M^{BM}(X/S)$ to $M^{BM}(Y/S)$ and denote $\alpha'=\epsilon_{X,Y}(\alpha)$. Then the following diagrams are commutative:
\begin{enumerate}
\item
\label{lemma_epsilon_1}
$$
  \xymatrix{
    M^{BM}(XY/B) \ar[r]^-{\alpha'} \ar[d]_-{R(\alpha)} 
& \mathbbold{1}_B(-d_Y)[-2d_Y] \\
    M^{BM}(YY/B) \ar[r]^-{\delta_{Y/S}^*} & M^{BM}(Y/B) \ar[u]_-{G(Y/B)}.
  }
$$
\item
\label{lemma_epsilon_2}
$$
  \xymatrix{
    M^{BM}(X/S) \ar[r]^-{\alpha} \ar[d]_-{p^*} 
&  M^{BM}(Y/S)\\
    M^{BM}(XY/S) \ar[r]^-{R(\delta_{Y/B})} & M^{BM}(XY\times Y/S)(d_Y)[2d_Y], \ar[u]_-{R(\alpha')}
  }
$$
where $p:XY\to X$ is the canonical projection.
\end{enumerate}
\end{lemma}

\proof
By definition, $\alpha'$ is the composition 
\begin{align*}
        &M^{BM}(XY/B)\overset{BC}{\simeq} h_!g^*M^{BM}(X/S)\overset{\alpha}{\to} h_!g^*M^{BM}(Y/S)\\
        &\xrightarrow{ad'_{(g^*,g_*)}}M^{BM}(Y/B)\xrightarrow{G(Y/B)}\mathbbold{1}_B(-d_Y)[-2d_Y],
\end{align*}
and the first diagram follows from Lemma~\ref{compat_BC_diag}. The second diagram is similar.
\endproof

\begin{proposition}
\label{composition_law}
Let $X$, $Y$ and $Z$ be three projective $S$-schemes that are smooth over $B$, and let $\alpha$ (respectively $\beta$) be a map from $M^{BM}(X/S)$ to $M^{BM}(Y/S)$ (respectively from $M^{BM}(Y/S)$ to $M^{BM}(Z/S)$). We have the following cartesian diagram:
$$
  \xymatrix{
    XYZ \ar[r]^-{} \ar[d]^-{} & XY\times YZ \ar[d]^-{f}\\
    Y \ar[r]^-{\delta_{Y/B}} & Y\times Y.
  }
$$
Denote by $p:XYZ\to XZ$ the canonical projection.
Then we have the following equality:
$$
\epsilon_{X,Z}(\beta\circ\alpha)=p^*\circ R_f(\delta_{Y/B})(\epsilon_{X,Y}(\alpha)\otimes\epsilon_{Y,Z}(\beta))
$$
\end{proposition}

\proof
Denote respectively $\alpha'=\epsilon_{X,Y}(\alpha)$ and $\beta'=\epsilon_{Y,Z}(\beta)$. We need to show that the following diagram is commutative:
$$
  \xymatrix{
    M^{BM}(XZ/B) \ar[r]^-{p^*} \ar[rd]_-{R(\alpha)} & M^{BM}(XYZ/B) \ar@{}[d]|-{(1)} \ar[r]^-{R(\delta_{Y/B})} & M^{BM}(XY\times YZ/B)(d_Y)[2d_Y]  \ar[d]^-{\alpha'\otimes\beta'} \ar[ld]_-{\alpha'\otimes\mathbbold{1}} \\
   & M^{BM}(YZ/B) \ar@{}[rd]|-{(3)} \ar[r]^-{\beta'} \ar[d]_-{R(\beta)} & \mathbbold{1}_B(-d_Z)[-2d_Z] \ar@{}[]+L+<-6pt,+20pt>|-{(2)} \\
& M^{BM}(ZZ/B) \ar[r]^-{\delta^*_Z} & M^{BM}(Z/B). \ar[u]_-{G(Z/B)}
  }
$$
The triangle $(2)$ is commutative by Lemma~\ref{kunneth}, and the square $(3)$ follows from Lemma~\ref{lemma_epsilon} (\ref{lemma_epsilon_1}). For the square $(1)$, we are reduced to the following diagram:
$$
  \xymatrix{
    M^{BM}(XZ/S) \ar[r]^-{R(\alpha)} \ar[d]_-{p^*} 
&  M^{BM}(YZ/S)\\
    M^{BM}(XYZ/S) \ar[r]^-{R(\delta_{Y/B})} & M^{BM}(XY\times YZ/S)(d_Y)[2d_Y], \ar[u]_-{R(\alpha')=\alpha'\otimes\mathbbold{1}}
  }
$$
which is the base change of the diagram in Lemma~\ref{lemma_epsilon} (\ref{lemma_epsilon_2}) by the proper morphism $Z\to S$. Therefore the result follows from the transitivity of the $R$-operation (Lemma~\ref{compat_r}) and the compatibility of $R$-operation with proper functoriality (Lemma~\ref{lemma_monoidal_BC} (\ref{lemma_monoidal_BC1})).
\endproof

\section{Borel-Moore homology and Chow groups}
\label{Chapter3}
\noindent
In this chapter, $k$ is a perfect base field, and all $k$-schemes are supposed to be seperated of finite type. According to this convention, a $k$-scheme is regular if and only if it is smooth over $k$, and every $k$-scheme contains a smooth open subscheme. For any abelian group $A$, denote by $A_\mathbb{Q}$ the $\mathbb{Q}$-vector space $A\otimes_\mathbb{Z}\mathbb{Q}$.
\subsection{The niveau spectral sequence}\ 

Let $S$ be a scheme and $X$ be a separated $S$-scheme of finite type. Recall that (Definition~\ref{BM_homology}) the Borel-Moore homology of $X$ relative to $S$ is defined as
$$
H^{BM}_{m,n}(X/S):=Hom_{DM_{{\cyrrm{B}},c}(S)}(M^{BM}(X/S),\mathbbold{1}_S(-n)[-m]).
$$
The aim of this section is to identify Borel-Moore homology with Chow groups, using the classical tool of coniveau spectral sequence. We follow the treatment in \cite{CTHK} and \cite{Deg1}.
Note that all the functorialities of Borel-Moore motives established in Chapter~\ref{chapter2}, including
\begin{itemize}
\item proper and smooth functorialities
\item localization distinguished triangle
\item K\"unneth formula
\item Gysin and refined Gysin morphisms
\item specialization maps
\end{itemize}
become naturally functorialities for the Borel-Moore homology, in a contravariant way. 
If $X$ is a $k$-scheme, denote by $M^{BM}(X)=M^{BM}(X/k)$ and $H^{BM}_{p,q}(X)=H^{BM}_{p,q}(X/k)$. 
The localization distinguished triangle (Lemma~\ref{bm_basic} (\ref{localization_leq})) leads to a long exact sequence of localization on Borel-Moore homology:
\begin{corollary}
\label{leq_BM}
Let $X$ be a $k$-scheme and $i:Z\to X$ be a closed immersion with complementary open immersion $j:X-Z\to X$. There is a long exact sequence
$$
\cdots\to H^{BM}_{p,q}(Z)\overset{i_*}{\to} H^{BM}_{p,q}(X)\overset{j^*}{\to} H^{BM}_{p,q}(X-Z)\overset{\partial_{X,Z}}{\to} H^{BM}_{p-1,q}(Z)\to\cdots.
$$

\end{corollary}


\begin{definition}
Let $X$ be a scheme. A \textbf{tower} on $X$ is an increasing chain $\bar{Z}=(Z_i)_i$ of closed subschemes of $X$ of the form
$$
Z_{-1}=\emptyset\subset Z_0\subset Z_1\subset\cdots
$$
such that $dim\ Z_i\leqslant i$ and $Z_i=Z_{i+1}$ for sufficiently large $i$. We denote by $Z_\infty$ the increasing limit (=union) of the chain $\bar{Z}$.
\end{definition}

Now we construct the niveau spectral sequence, following the same steps as \cite[1.2]{CTHK} and \cite[Definition 1.6]{Deg1}. Fix a $k$-scheme $X$ and let $\bar{Z}$ be a tower on $X$.
Applying Lemma~\ref{leq_BM} to the closed immersion $Z_{p-1}\to Z_p$, we have a long exact sequence
$$
\cdots\to H^{BM}_{m,n}(Z_{p-1})\to H^{BM}_{m,n}(Z_p)\to H^{BM}_{m,n}(Z_p-Z_{p-1})\to H^{BM}_{m-1,n}(Z_{p-1})\to\cdots,
$$
which provides an exact couple $(D,E)$, with $D_{p,q}(\bar{Z},i)=H^{BM}_{p+q,i}(Z_p)$ and  $E_{p,q}(\bar{Z},i)=H^{BM}_{p+q,i}(Z_p-Z_{p-1})$.
\footnote{The index $i$ here means the $i$-th graded part of the total term.}
Such an exact couple gives rise to a spectral sequence (of homological type):
$$
E^1_{p,q}(\bar{Z},i)=H^{BM}_{p+q,i}(Z_p-Z_{p-1})
$$
converging to $H^{BM}_{p+q,i}(Z_\infty)$, with respect to the increasing filtration $F_p=Im(H^{BM}_{p+q,i}(Z_p)\to H^{BM}_{p+q,i}(Z_{\infty}))$. (i.e. $E^\infty_{p,q}(\bar{Z},i)=F_p/F_{p-1}$).

Now by taking the limit of all towers $\bar{Z}$, we obtain an exact couple, which we denote by $(D_{p,q}(X,i),E_{p,q}(X,i))$,
and therefore a spectral sequence of the form
$$
E^1_{p,q}(X,i)=\underset{\longrightarrow}{lim}\ H^{BM}_{p+q,i}(Z_p-Z_{p-1})\Rightarrow H^{BM}_{p+q,i}(X),
$$
called respectively the \textbf{niveau exact couple} and the \textbf{niveau spectral sequence} of $X$.

\begin{definition}

Let $x\in X$ be a topological point. Define
$$
\hat{H}^{BM}_{n,i}(x):=\underset{\longrightarrow}{lim}\ H^{BM}_{n,i}(\bar{x}\cap U)
$$
where the limit is taken over Zariski open neighborhoods $U$ of $x$ in $X$.
\footnote{In view of the theory of generic motives, $\hat{H}^{BM}_{n,i}(x)$ is nothing else but the Borel-Moore homology of the generic motive associated to the residue field of $x$, which only depends on the residue field, see \cite{Deg1}.}

\end{definition}
By Remark~\ref{BM_union} (\ref{BM_union_1}), the $E^1$ terms can be identified as 
$$
E^1_{p,q}(X,i)=\underset{\longrightarrow}{lim}\ H^{BM}_{p+q,i}(Z_p-Z_{p-1})=\underset{x\in X_{(p)}}{\oplus} \hat{H}^{BM}_{p+q,i}(x)
$$

The following lemma results from Lemma~\ref{lemma_voe} and Lemma~\ref{ident_BM_mot}, which shows that the spectral sequence is bounded:
\begin{lemma}
If $x$ is a $d$-dimensional point of $X$, then $\hat{H}^{BM}_{p,q}(x)=0$ whenever $p<d+q$ or $q>d$.
\end{lemma}

\begin{corollary}
The term $E^1_{p,q}(X,i)$ is zero if $q<i$ or $p<i$.
\end{corollary}

By the previous corollary, the group $H^{BM}_{2n,n}(X)$ can be identified as the cokernel of the map
$$
d^1_{n+1,n}:E^1_{n+1,n}(X,n)\to E^1_{n,n}(X,n).
$$

The following lemma, which will be proved in Section~\ref{proof_elem_BM}, 
describes the local behavior of basic operations on Borel-Moore homology:
\begin{lemma}
\label{elem_BM}
For any irreducible smooth $k$-scheme $X$ of dimension $d$, there are isomorphisms $H^{BM}_{2d,d}(X)\simeq\mathbb{Q}$ and $H^{BM}_{2d-1,d-1}(X)\simeq \mathcal{O}^*(X)_\mathbb{Q}$. Consequently, for any $k$-scheme $X$ and any $d$-dimensional point $x$ of $X$, there are isomorphisms
$\hat{H}^{BM}_{2d,d}(x)\simeq\mathbb{Q}$ and $\hat{H}^{BM}_{2d-1,d-1}(x)\simeq\kappa(x)^*_\mathbb{Q}$, where $\kappa(x)$ is the residue field of $x$. In addition, these isomorphisms satisfy the following functorial properties:
\begin{enumerate}
\item 
\label{elem_BM_proper}
For any finite surjective morphism $p:X\to Y$ between irreducible smooth $k$-schemes both of dimension $d$, the proper functoriality map $p_*:H^{BM}_{2d,d}(X)\to H^{BM}_{2d,d}(Y)$ on Borel-Moore homology is the multiplication by the degree of field extension $[k(X):k(Y)]$ as a map $\mathbb{Q}\to\mathbb{Q}$;

\item 
\label{elem_BM_smooth}
For any smooth morphism $f:X\to Y$ between irreducible smooth $k$-schemes, of dimension respectively $d$ and $d'$, the smooth functoriality map $f^*:H^{BM}_{2d',d'}(Y)\to H^{BM}_{2d,d}(X)$ is the identity map $\mathbb{Q}\to\mathbb{Q}$;

\item 
\label{elem_BM_boundary}
For any closed immersion $i:Z\to X$ of codimension $1$ between irreducible smooth $k$-schemes ($Z$ is then a divisor of $X$), with $x$ the generic point of $X$ and $d$ the dimension of $X$, the following diagram is commutative:
$$
  \xymatrix{
    H^{BM}_{2d-1,d-1}(X-Z)\simeq \mathcal{O}^*(X-Z)_\mathbb{Q} \ar[rr]^-{\partial_{X,Z}} \ar[rd] & & H^{BM}_{2d-2,d-1}(Z)\simeq\mathbb{Q},\\
       &  \hat{H}^{BM}_{2d-1,d-1}(x)\simeq\kappa(x)^*_\mathbb{Q} \ar[ru]^-{ord_Z} &
  }
$$
where $ord_Z:\kappa(x)^*_\mathbb{Q}\to\mathbb{Q}$ is the order of vanishing along $Z$, extended to rational coefficients (\cite[1.2]{Ful}).

\item 
\label{comp_prod} 
Let $X$ and $Y$ be irreducible smooth $k$-schemes of dimension respectively $d$ and $d'$. Then the tensor product map $H^{BM}_{2d,d}(X)\otimes H^{BM}_{2d',d'}(Y)\to H^{BM}_{2d+2d',d+d'}(X\times Y)$ defined by the K\"unneth isomorphism (Lemma~\ref{kunneth}) sends the class $1\otimes 1$ to the class $\sum_I1$, where elements of $I$ are irreducible components of $X\times Y$.
\end{enumerate}
\end{lemma}


Therefore we have 
$$
E^1_{n+1,n}(X,n)=\underset{x\in X_{(n+1)}}{\oplus} \hat{H}^{BM}_{2n+1,n}(x)\simeq\underset{x\in X_{(n+1)}}{\oplus}(\kappa(x))^*_\mathbb{Q}
$$
and 
$$
E^1_{n,n}(X,n)=\underset{x\in X_{(n)}}{\oplus} \hat{H}^{BM}_{2n,n}(x)\simeq\underset{x\in X_{(n)}}{\oplus}\mathbb{Q}=\mathcal{Z}_n(X)_{\mathbb{Q}},
$$
where $\mathcal{Z}_n(X)$ is the group of $n$-dimensional algebraic cycles of $X$ with rational coefficients. The next objective is to study the map $d^1_{n+1,n}$.

\begin{definition}

Let $X$ be a scheme and $z\in X_{(n+1)}$ and $y\in X_{(n)}$ be two points. Suppose that $y$ is a specialization of $z$. Let $Z$ be the reduced closure of $z$ in $X$ and $\tilde{Z}\overset{f}{\to} Z$ its normalization. Let $z'=f^{-1}(z)$. 
Let $Z_y$ be the reduced closure of $y$ in $X$ and for each $t\in f^{-1}(y)$, let $Z_t$ be the reduced closure of $t$ in $\tilde{Z}$. Then $f$ induces a finite morphism $f_t:Z_t\to Z_y$.
We define the map
$$
\phi_{t*}:\hat{H}^{BM}_{p,q}(t)\to \hat{H}^{BM}_{p,q}(y)
$$
to be the limit of the maps 
$$
f_{t*|_{U}}: H^{BM}_{p,q}(Z_t\times_XU)\to H^{BM}_{p,q}(Z_y\times_XU),
$$
where $U$ runs through all Zariski open neighborhoods of $y$ in $X$, and the map
$$
\partial_t:\hat{H}^{BM}_{p+1,q}(z')\to \hat{H}^{BM}_{p,q}(t)
$$
to be the limit of the maps 
$$
\partial_{U,Z_t\times_XU}: H^{BM}_{p+1,q}(U-Z_t\times_XU)\to H^{BM}_{p,q}(Z_t\times_XU),
$$
where $U$ runs through all Zariski open neighborhoods of $t$ in $\tilde{Z}$. 

Composing these two maps, we obtain a map 
\footnote{The isomorphism is constructed in the following way: there is a natural map $\hat{H}^{BM}_{p,q}(z')\to\hat{H}^{BM}_{p,q}(z)$ defined in a similar way as the map $\phi^{*}_t$ above, and it is an isomorphism because $f$ induces an isomorphism between the residue fields of $z$ and $z'$, therefore induces an isomorphism between the associated pro-schemes.}
$$
res_t:\hat{H}^{BM}_{p+1,q}(z)\simeq
\hat{H}^{BM}_{p+1,q}(z')\overset{\phi_{t*}\circ\partial_{t}}{\to} \hat{H}^{BM}_{p,q}(y).
$$

We define the \textbf{residue map} as 
$$
res=\sum_{t\in f^{-1}(y)}res_t: \hat{H}^{BM}_{p+1,q}(z)\to \hat{H}^{BM}_{p,q}(y).
$$

If $y$ is not a specialization of $z$, we put $res=0$ as a map $\hat{H}^{BM}_{p+1,q}(z)\to \hat{H}^{BM}_{p,q}(y)$.
\end{definition}

\begin{lemma}
\label{diff_residue}
Let $z\in X_{(n+1)}$ and $y\in X_{(n)}$ be two points. Then 
\begin{enumerate}
\item We have a commutative diagram
$$
  \xymatrix{
    E^1_{n+1,n}(X,n) \ar[r]^-{d^1_{n+1,n}} \ar[d] & E^1_{n,n}(X,n)  \\
    \hat{H}^{BM}_{2n+1,n}(z) \ar[r]^-{res} & \hat{H}^{BM}_{2n,n}(y). \ar[u]
  }
$$
\item If $y$ is a specialization of $z$, then we have a commutative diagram
$$
  \xymatrix{
    \hat{H}^{BM}_{2n+1,n}(z) \ar[r]^-{res}  \ar[d]_-{\wr} & \hat{H}^{BM}_{2n,n}(y) \ar[d]^-{\wr} \\
   \kappa(z)^*_\mathbb{Q} \ar[r]^-{ord_y} & \mathbb{Q},
  }
$$
where both vertical isomorphisms are deduced from Lemma~\ref{elem_BM}.
\end{enumerate}
\end{lemma}
\proof

The proof is adapted from \cite[Proposition 1.16]{Deg1}. The differential map $d^1_{n+1,n}$ in the spectral sequence is the limit of maps
$$
H^{BM}_{2n+1,n}(Z_{n+1}-Z_n)\xrightarrow{\partial_{Z_{n+1},Z_n}}H^{BM}_{2n,n}(Z_n)\to H^{BM}_{2n,n}(Z_n-Z_{n-1}).
$$
Denote $Z=Z_{n+1}-Z_{n-1}$ and $Y=Z_n-Z_{n-1}$, and we know according to Lemma~\ref{comp_diff} that the previous map is equal to the map
$$
H^{BM}_{2n+1,n}(Z-Y)\xrightarrow{\partial_{Z,Y}}H^{BM}_{2n,n}(Y).
$$
As one takes the limit of towers, one may enlarge $Z_n$ and $Z_{n-1}$, and thus one may suppose that $Z-Y$ and $Y$ are regular schemes. We suppose $z\in Z$ and $y\in Y$. Denote by $Z_z$ (respectively $Y_y$) the irreducible component of $Z$ (respectively $Y$) containing $z$ (respectively $y$). 
As the scheme $Y$ is regular, the difference $Y'_y=Y-Y_y$ is the union of other connected components, thus is a closed subscheme of $Y$. Applying Lemma~\ref{comp_diff} again, we know that the composite map
$$
H^{BM}_{2n+1,n}(Z-Y)\xrightarrow{\partial_{Z,Y}}H^{BM}_{2n,n}(Y)\to H^{BM}_{2n,n}(Y_y).
$$
is equal to the map
$$
H^{BM}_{2n+1,n}(Z-Y)\xrightarrow{\partial_{Z-Y'_y,Y_y}}H^{BM}_{2n,n}(Y_y).
$$
We now study the limit of this map.

Suppose first that $y$ is not a specialization of $z$, then by enlarging $Z_{n-1}$, one may suppose that the scheme-theoretic intersection $Y_y\cap Z_z$ is empty.  Therefore we have a cartesian square of closed immersions
$$
  \xymatrix{
    \emptyset \ar[r] \ar[d] & Y_y \ar[d] \\
    Z_z-Y\cap Z_z \ar[r] & Z-(Y-Y_y).
  }
$$
According to Lemma~\ref{cartesian_diff_comp} $(1)$, the map $\partial_{Z-(Y-Y_y),Y_y}$ is zero, as the proposition asserts.

Now suppose that $y$ is a specialization of $z$, and we may suppose that the scheme $Z$ is irreducible with generic point $z$. Let $p:\tilde{Z}\to Z$ be the normalization map, and we know that the scheme $\tilde{Z}$ is regular in codimension $1$, so we may suppose that $\tilde{Z}$ is regular. Denote by $\tilde{Y}$ (respectively $\tilde{Y}_y$ and $\tilde{Y}'_y$) the fiber product $Y\times_Z \tilde{Z}$ (respectively $Y_y\times_Z \tilde{Z}$ and $Y'_y\times_Z \tilde{Z}$). We may suppose that $\tilde{Y}_y$ is regular and the two schemes $\tilde{Y}_y$ and $\tilde{Y}'_y$ are disjoint. In addition, we may suppose that $\tilde{Y}_y$ is of pure dimension $n$, so that all of its components dominates $Y_y$. Therefore there is a cartesian square
$$
  \xymatrix{
    \tilde{Y}_y \ar[r] \ar[d] & \tilde{Z}-\tilde{Y}'_y \ar[d] \\
    Y_y \ar[r] & Z-Y'_y.
  }
$$
According to Lemma~\ref{cartesian_diff_comp} $(1)$ there is a commutative square
$$
  \xymatrix{
    H^{BM}_{2n+1,n}(\tilde{Z}-\tilde{Y}) \ar[r]^-{\partial_{\tilde{Z}-\tilde{Y}'_y,\tilde{Y}_y}} \ar[d]_-{p_*} & H^{BM}_{2n,n}(\tilde{Y}_y) \ar[d]^-{p_*} \\
    H^{BM}_{2n+1,n}(Z-Y) \ar[r]^-{\partial_{Z-Y'_y,Y_y}}& H^{BM}_{2n,n}(Y_y).
  }
$$
The vertical map in the left is an isomorphism when we take the limit, so we are reduced to study the map $\partial_{\tilde{Z}-\tilde{Y}'_y,\tilde{Y}_y}$ in the upper row. The scheme $\tilde{Y}_y$ is not necessarily connected; however, the set of its connected components are in bijection with the set $f^{-1}(y)$, and for any $t\in f^{-1}(y)$ we denote by $\tilde{Y}_t$ the corresponding component, and $\tilde{Y}'_t=\tilde{Y}-\tilde{Y}_t$. For every $t\in f^{-1}(y)$ there is a finite surjective map $p_t:\tilde{Y}_t\to Y_y$ between integral schemes. Then by Lemma~\ref{comp_diff} and Remark~\ref{BM_union} (\ref{BM_union_1}), the following diagram is commutative:
$$
  \xymatrix{
    H^{BM}_{2n+1,n}(\tilde{Z}-\tilde{Y}) \ar[r]^-{\partial_{\tilde{Z}-\tilde{Y}'_y,\tilde{Y}_y}} \ar[rd]_-{\underset{t}{\oplus}(\partial_{\tilde{Z}-\tilde{Y}'_t,\tilde{Y}_t})} & H^{BM}_{2n,n}(\tilde{Y}_y) \ar[r]^-{p_*} & H^{BM}_{2n,n}(Y_y), \\
     & \underset{t\in f^{-1}(y)}{\oplus}H^{BM}_{2n,n}(\tilde{Y}_t) \ar[u]_-{\wr} \ar[ru]_-{\sum_t p_{t*}} &
  }
$$
which proves the first part of the lemma. For the second part, it follows from Lemma~\ref{elem_BM} that the following diagram is commutative:
$$
  \xymatrix{
    H^{BM}_{2n+1,n}(\tilde{Z}-\tilde{Y}) \ar[rr]^-{\partial_{\tilde{Z}-\tilde{Y}'_t,\tilde{Y}_t}} \ar[rd] & & H^{BM}_{2n,n}(\tilde{Y}_t), \\ 
     &  \hat{H}^{BM}_{2n+1,n}(z)=\kappa(z)^*_\mathbb{Q} \ar[ru]_-{ord_{\tilde{Y}_t}} & 
  }
$$
and the result is straightforward.

\endproof

\begin{corollary}
\label{BM=Chow}
For any separated $k$-scheme of finite type $X$, there is a canonical isomorphism
$$
H^{BM}_{2n,n}(X)\simeq CH_n(X)_{\mathbb{Q}}.
$$
\end{corollary}
\proof

Lemma~\ref{diff_residue} 
shows that the map $d^1_{n+1,n}:E^1_{n+1,n}(X,n)\to E^1_{n,n}(X,n)$ is the divisor class map 
$$
div:\underset{x\in X_{(n+1)}}{\oplus}(\kappa(x))^*_\mathbb{Q}\to\underset{x\in X_{(n+1)}}{\oplus}\mathbb{Q},
$$
whose cokernel is nothing else but the Chow group.

\endproof
In the next section we are going to study functorial properties of this isomorphism.

\subsection{Operations on Chow groups}\

We now show that the isomorphism in Corollary~\ref{BM=Chow} is compatible with proper and smooth functorialities and refined Gysin maps on Chow groups. 
Start with the special case of open immersions:
\begin{lemma}
\label{open_immersion_identity}
Let $X$ be an irreducible $n$-dimensional $k$-scheme and $U$ a nonempty open subscheme of $X$. Then the pull-back map $H^{BM}_{2n,n}(X)\to H^{BM}_{2n,n}(U)$ is the identity map $\mathbb{Q}\to\mathbb{Q}$.
\end{lemma}
\proof
We know that $H^{BM}_{2n,n}(X)\simeq E^1_{n,n}(X,n)$ and $H^{BM}_{2n,n}(U)\simeq E^1_{n,n}(U,n)$. It is not difficult to check that the pull-back map agrees with the map $E^1_{n,n}(X,n)\to H^{BM}_{2n,n}(U)$ defined by taking limits of the maps $H^{BM}_{2n,n}(X-Z_{n-1})\to H^{BM}_{2n,n}((X-Z_{n-1})\cap U)$. Taking $Z_{n-1}$ large enough we may suppose that $X-Z_{n-1}$ is regular, and the result follows from Lemma~\ref{elem_BM} (\ref{elem_BM_smooth}).
\endproof

\begin{proposition}
\label{compat_proper_pf}
For any proper morphism $p:X\to Y$ between separated $k$-schemes of finite type, the diagram
$$
  \xymatrix{
    H^{BM}_{2n,n}(X)\ar[r]^-{p_*} \ar[d]^-{\wr} & H^{BM}_{2n,n}(Y) \ar[d]^-{\wr} \\
    CH_n(X)_\mathbb{Q} \ar[r]^-{p_*}  & CH_n(Y)_\mathbb{Q}
  }
$$
is commutative.
\end{proposition}

\proof

The push-forward on Chow groups $CH_n(X)\overset{p_*}{\to}CH_n(Y)$ is induced by the push-forward on the group of cycles of dimension $n$: $\mathcal{Z}_n(X)\overset{p_*}{\to}\mathcal{Z}_n(y)$. In view of Lemma~\ref{diff_residue}, we construct a push-forward map on the niveau spectral sequence, and we show that the map induced on the group of cycles is the one we need.

Let $Z_k$ be a closed subscheme of $X$ of dimension less than or equal to $k$. As $p:X\to Y$ is proper, the image $Z_k'=p(Z_k)$ is a closed subscheme of $Y$ of dimension at most $k$, and the restriction of $p$ is a proper surjective morphism $p_1:Z_k\to Z_k'$. Let $Z_{k-1}'$ be a closed subscheme of $Z_k'$ which does not contain the image of any $k$-dimensional point of $Z_k$ by $p_1$, and let $Z_{k-1}$ be the fiber product $Z_{k-1}'\times_{Z_k'}Z_k$. Then both $Z_{k-1}$ and $Z_{k-1}'$ are of dimension less than or equal to $k-1$, and we have a commutative diagram
$$
  \xymatrix{
    Z_{k-1} \ar[r] \ar[d]^-{} & Z_k \ar[d]^-{p_1} & Z_k-Z_{k-1} \ar[l]_-{j} \ar[d]^-{q_1} \\
    Z_{k-1}' \ar[r] & Z_k' & Z_k'-Z_{k-1}', \ar[l]_-{j'}
  }
$$
where both squares are cartesian.
We can complete the couple $(Z_k,Z_{k-1})$ (respectively the couple $(Z_k',Z_{k-1}')$) into a tower $\bar{Z}$ (respectively $\bar{Z'}$), such that the previous relation is satisfied for all $k$. Then for all $r$ and $s$, the maps 
$$
H^{BM}_{r,s}(Z_k) \overset{p_{1*}}{\to}H^{BM}_{r,s}(Z_k')
$$
$$
\textrm{and }H^{BM}_{r,s}(Z_k-Z_{k-1}) \overset{q_{1*}}{\to}H^{BM}_{r,s}(Z_k'-Z_{k-1}')
$$
define respectively maps $D^1_{r,s}(\bar{Z},i)\to D^1_{r,s}(\bar{Z'},i)$ and $E^1_{r,s}(\bar{Z},i)\to E^1_{r,s}(\bar{Z'},i)$. Using Lemma~\ref{BC_proper_smooth}, we can show that these maps define a morphism of exact couples. On the other hand, we can show that the exact couple obtained by taking limit of towers of the form $\bar{Z}$ is isomorphic to the niveau exact couple of $X$. Therefore there is a well-defined map from the niveau exact couple of $X$ to that of $Y$.

Now we focus in particular on the induced maps on the groups $H^{BM}_{2n,n}$. Without loss of generality, we may suppose that $Z_n$ is irreducible, and therefore $Z_n'=p(Z_n)$ is irreducible as well.
Then by Lemma~\ref{BC_proper_smooth} there is a commutative diagram
$$
  \xymatrix{
H^{BM}_{2n,n}(Z_n) \ar[r]^-{j^*} \ar[d]^-{p_{1*}} &  H^{BM}_{2n,n}(Z_n-Z_{n-1})  \ar[d]^-{q_{1*}}\\ 
H^{BM}_{2n,n}(Z_n')  \ar[r]^-{j'^*} &  H^{BM}_{2n,n}(Z_n'-Z_{n-1}'),  
  }
$$
where by Lemma~\ref{open_immersion_identity}, the two horizontal maps are both isomorphisms.
In consequence, we have the following commutative diagram
$$
  \xymatrix{
    H^{BM}_{2n,n}(Z_n-Z_{n-1}) \ar[d]^-{q_{1*}} \ar[r]^-{\sim} & H^{BM}_{2n,n}(Z_n) \ar[d]^-{p_{1*}} \ar[r] & H^{BM}_{2n,n}(X)  \ar[d]^-{p_*}\\
   H^{BM}_{2n,n}(Z_n'-Z_{n-1}')  \ar[r]^-{\sim} & H^{BM}_{2n,n}(Z_n')  \ar[r] & H^{BM}_{2n,n}(Y).
  }
$$
By taking limit on the first column, we obtain a map $p_*:E^1_{n,n}(X,n)\to E^1_{n,n}(Y,n)$, which fits in a commutative diagram
$$
  \xymatrix{
    E^1_{n,n}(X,n) \ar[r]^-{p_{*}} \ar[d] & E^1_{n,n}(Y,n) \ar[d] \\
    H^{BM}_{2n,n}(X) \ar[r]^-{p_*} & H^{BM}_{2n,n}(Y).
  }
$$

It remains to compute the map $p_*$ on the first row, which is the limit of the maps $q_{1*}:H^{BM}_{2n,n}(Z_n-Z_{n-1})\to H^{BM}_{2n,n}(Z'_n-Z'_{n-1})$ . If the dimension of $Z'_n=p(Z_n)$ is less than $n$, then $H^{BM}_{2n,n}(Z_n'-Z_{n-1}')=0$ and the map is zero.

Suppose now $Z'_n$ is $n$-dimensional, then the map $q_{1*}$ is finite and surjective. 
By enlarging $Z'_{n-1}$, one may suppose that both $Z_n-Z_{n-1}$ and $Z'_n-Z'_{n-1}$ are regular. Then by Lemma~\ref{elem_BM} (\ref{elem_BM_proper}) the map $q_{1*}:H^{BM}_{2n,n}(Z_n-Z_{n-1})\to H^{BM}_{2n,n}(Z'_n-Z'_{n-1})$ is the multiplication by the degree $[k(Z'_n-Z'_{n-1}):k(Z_n-Z_{n-1})]$ of the extension between the relevant function fields as a map $\mathbb{Q}\to\mathbb{Q}$.

In the general case on deduces from the discussion above that the map $E^1_{n,n}(X,n)\to E^1_{n,n}(Y,n)$ is the push-forward map on cycles, which completes the proof.
\endproof

\begin{proposition}
\label{comp_smooth_pullback}
For any smooth morphism $f:Y\to X$ of relative dimension $d$ between separated $k$-schemes of finite type, the diagram
$$
  \xymatrix{
    H^{BM}_{2n,n}(X) \ar[r]^-{f^*} & H^{BM}_{2n+2d,n+d}(Y) \\
    CH_n(X)_\mathbb{Q} \ar[u]^-{\wr} \ar[r]^-{f^*} & CH_{n+d}(Y)_\mathbb{Q} \ar[u]_-{\wr} 
  }
$$
is commutative, where the map on the second line is the usual pull-back on Chow groups. 
\end{proposition}
\proof

In the same way as for the previous proposition, we start with the construction of a pull-back map on the niveau spectral sequence. Let $Z_n$ be an $n$-dimensional irreducible closed subscheme of $X$ and $Z_{n-1}$ be a closed subscheme of $Z_n$ of pure dimension $(n-1)$. As the morphism $f:X\to Y$ is flat of relative dimension $d$, the scheme-theoretic preimages $f^{-1}(Z_n)$ and $f^{-1}(Z_{n-1})$ are both closed subschemes of $X$ respectively of pure dimension $n+d$ and $n+d-1$. Therefore in the diagram
$$
  \xymatrix{
    f^{-1}(Z_n)-f^{-1}(Z_{n-1})\ar@{}[rd]|{(1)} \ar[r] \ar[d] & f^{-1}(Z_n) \ar@{}[rd]|{(2)} \ar[r] \ar[d] & Y \ar[d]^-{f} \\
    Z_n-Z_{n-1} \ar[r] & Z_n \ar[r] & X,
  }
$$
both squares $(1)$ and $(2)$ are cartesian. Applying functoriality to the square $(1)$ and Lemma~\ref{BC_proper_smooth} to the square $(2)$, we obtain the following commutative diagram
$$
  \xymatrix{
    H^{BM}_{2n+2d,n+d}(f^{-1}(Z_n)-f^{-1}(Z_{n-1})) & H^{BM}_{2n+2d,n+d}(f^{-1}(Z_n)) \ar[l]_-{(*)} \ar[r]  & H^{BM}_{2n+2d,n+d}(Y)  \\
    H^{BM}_{2n,n}(Z_n-Z_{n-1}) \ar[u]_-{(3)} & H^{BM}_{2n,n}(Z_n) \ar[l]_-{(*)} \ar[r] \ar[u] & H^{BM}_{2n,n}(X). \ar[u]^-{f^*}
  }
$$
By Lemma~\ref{open_immersion_identity}, the two maps marked with $(*)$ are isomorphisms, and therefore by 
taking the limit of the map of the form $(3)$ with respect to the tower $\bar{Z}$, we obtain a map between the $E^1$-terms of the niveau spectral sequence, with a commutative diagram
$$
  \xymatrix{
    E^1_{n+d,n+d}(Y,n+d) \ar[r]  & H^{BM}_{2n+2d,n+d}(Y)  \\
    E^1_{n,n}(X,n) \ar[u] \ar[r] \ar[u] & H^{BM}_{2n,n}(X). \ar[u]^-{f^*}
  }
$$
Then we only need to study the map $(3)$ for $Z_n$ and $Z_{n-1}$ large enough. One may suppose that both schemes $f^{-1}(Z_n)-f^{-1}(Z_{n-1})$ and $Z_n-Z_{n-1}$ are regular. Then by Lemma~\ref{elem_BM} (\ref{elem_BM_smooth}), the map $E^1_{n,n}(X,n)\to E^1_{n+d,n+d}(Y,n+d)$ is the diagonal embedding $\mathbb{Z}\to\mathbb{Z}^{b_0(f^{-1}(Z_n))}$.
In other words, for any $n$-dimensional point $x$ of $X$, the map $f^*$ sends the class of $x$ in $H^{BM}_{2n,n}(X)$ to the cycle $\sum_{y\in f^{-1}(x)} y$ in $H^{BM}_{2n+2d,n+d}(Y)$. However, since the morphism $f$ is smooth, for any point $y\in f^{-1}(x)$, its Zariski closure $\bar{y}$ in $Y$ is smooth over $\bar{x}$, which is an integral scheme. Therefore the geometric multiplicity (\cite[1.5]{Ful}) of $\bar{y}$ in $Y$ is $1$, and the proposition follows since we have checked that $f^*$ is exactly the pull-back map on cycles.
\endproof

\begin{proposition}
\label{compat_sp}
Let $i:Z\to X$ be a closed immersion between separated $k$-schemes of finite type and $C$ be the corresponding normal cone of $Z$ in $X$. Then the following diagram is commutative
$$
  \xymatrix{
    H^{BM}_{2n,n}(X) \ar[r]^-{sp} \ar[d]_-{\wr} & H^{BM}_{2n,n}(C) \ar[d]^-{\wr} \\
    CH_n(X)_\mathbb{Q} \ar[r]^-{sp} & CH_{n}(C)_\mathbb{Q},
  }
$$
where the map $sp$ in the lower row is the specialization map on Chow groups (\cite[5.2]{Ful}).
\end{proposition}

\proof

First we need to construct a map between niveau spectral sequences. For any morphism $Y\to X$, denote by $C_Y$ the normal cone of $Y\times_X Z$ in $Y$. Then there is a corresponding specialization map $M^{BM}(C_Y)\to M^{BM}(Y)$ (Definition~\ref{def_specialization}).

Let $\bar{Z}=(Z_p)$ be a tower of closed subschemes of $X$. Let $C'_p$  be the fiber product $Z_{p-1}\times_{Z_p}C_{Z_p}$, which is a closed subscheme of $C_{Z_p}$. Then $C_{Z_{p-1}}$ is a closed subcone of $C'_p$, and $C_{Z_p-Z_{p-1}}$ is the complementary open subscheme of $C'_p$ in $C_{Z_p}$. Thus we get the commutative diagram:
$$
  \xymatrix{
    M^{BM}(Z_p-Z_{p-1}) \ar[r] \ar@{}[rd]|{(1)} & M^{BM}(Z_p) \ar[r] \ar@{}[rd]|{(2)} & M^{BM}(Z_{p-1}) \ar[r] & M^{BM}(Z_p-Z_{p-1})[1] \\
    M^{BM}(C_{Z_p-Z_{p-1}}) \ar[r] \ar[u] & M^{BM}(C_{Z_p}) \ar[r] \ar[u] & M^{BM}(C'_p) \ar[r]^-{(*)} \ar[u] & M^{BM}(C_{Z_p-Z_{p-1}})[1], \ar[u]
  }
$$
where both rows are localization distinguished triangles, and vertical maps are specialization maps. Indeed, by Lemma~\ref{uniqueness_h} we need only to show the commutativity of the two squares $(1)$ and $(2)$, which follows from Lemma~\ref{cartesian_diff_comp}. 

We know that if $Z_p$ is purely $d$-dimensional, then $C_{Z_p}$ is $d$-dimensional as well (\cite[B.6.6]{Ful}). Therefore in any case  the dimension of $C_{Z_p}$ is at most $p$, and that of $C'_p$ is at most $p$ as well. Let $T$ be the union of all $p$-dimensional components of $C'_p$ (possibly empty). Then $T$ is a union of connected components of $C_{Z_p}$ for dimension reason. 
Denote by $C"_p=C'_p-T$, whose dimension is at most $p-1$. 
Then by Lemma~\ref{compat_diff_two}, the map $(*)$ factors as $M^{BM}(C'_p)\to M^{BM}(C"_p)\to M^{BM}(C_{Z_p-Z_{p-1}})[1]$.
Again by Lemma~\ref{compat_diff_two}, there is a commutative triangle
$$
  \xymatrix{
    M^{BM}(Z_{p-1})& \\
    M^{BM}(C'_p) \ar[r]^-{i^*} \ar[u]^-{sp} & M^{BM}(C"_p), \ar[lu]_-{sp'}
  }
$$
where the map $sp'$ is defined in the same fashion as the specialization map. Finally we obtain a commutative diagram
$$
  \xymatrix{
    M^{BM}(Z_p-Z_{p-1}) \ar[r] & M^{BM}(Z_p) \ar[r] & M^{BM}(Z_{p-1}) \ar[r] & M^{BM}(Z_p-Z_{p-1})[1] \\
    M^{BM}(C_{Z_p-Z_{p-1}}) \ar[r] \ar[u] & M^{BM}(C_{Z_p}) \ar[r] \ar[u] & M^{BM}(C"_p) \ar[r] \ar[u] & M^{BM}(C_{Z_p-Z_{p-1}})[1], \ar[u]
  }
$$
with $C"_p$ of dimension at most $p-1$. The corresponding maps on cohomology defines a map from the niveau spectral sequence of $X$ to that of $C$. By Lemma~\ref{sp_identity} below, the induced map on $E^2$ is the specialization map on Chow groups, and the result follows.




\endproof

\begin{lemma}
\label{sp_identity}
Let $Z\to X$ be a closed immersion between irreducible regular schemes, where $X$ is $n$-dimensional, and let $C$ be the normal cone of $Z$ in $X$. Then the specialization map $H^{BM}_{2n,n}(X)\overset{sp}{\to}H^{BM}_{2n,n}(C)$ is the identity map $\mathbb{Q}\to\mathbb{Q}$.

\end{lemma}

\proof

By construction, the specialization map factors as $H^{BM}_{2n,n}(X)\to H^{BM}_{2n+1,n}(X\times\mathbb{G}_m)\overset{\partial_{D_ZX,C}}{\longrightarrow} H^{BM}_{2n,n}(C)$. 
The first map is the canonical inclusion $\mathbb{Q}\to\mathbb{Q}\oplus\mathcal{O}^*(X)_\mathbb{Q}$, which sends the class $1$ to $(1,0)$. By Lemma~\ref{elem_BM} (\ref{elem_BM_boundary}), the image of $(1,0)$ by $\partial_{D_ZX,C}$ is the vanishing order of the fundamental parameter of $\mathbb{G}_m$ along $C$, which equals to $1$, and the result follows.

\endproof

\begin{proposition}
\label{compat_refined_gysin}
Let $X'$ and $Y'$ be as in Definition~\ref{refined_gysin}. Then the following diagram is commutative
$$
  \xymatrix{
    H^{BM}_{2n,n}(Y') \ar[r]^-{R_f(i)} \ar[d]^-{\wr} & H^{BM}_{2n-2c,n-c}(X') \ar[d]^-{\wr} \\
    CH_n(Y')_\mathbb{Q} \ar[r]^-{R_f(i)} & CH_{n-c}(X')_\mathbb{Q},
  }
$$
where the map $R_f(i)$ in the lower row is the refined Gysin map on Chow groups (\cite[6.2]{Ful}).
\end{proposition}

\proof
Let $C'$ be the normal cone of $X'$ in $Y'$. By definition, the refined Gysin map on Chow groups is the composite map $CH_n(Y')\overset{sp}{\to} CH_n(C')\overset{i_*}{\to} CH_n(N')\overset{(\pi'^*)^{-1}}{\longrightarrow} CH_{n-c}(X')$ (since the map $N'\overset{\pi'}{\to} X'$ is a vector bundle, the pull-back map $\pi'^*$ between Chow groups is an isomorphism by \cite[Theorem 3.3(a)]{Ful}). Therefore it is sufficient to show that the following diagram is commutative:
$$
  \xymatrix{
    CH_n(Y')_\mathbb{Q} \ar[r]^-{sp} \ar[d]^-{\wr} \ar@{}[rd]|{(1)} & CH_n(C')_\mathbb{Q} \ar[r]^-{i_*} \ar[d]^-{\wr} \ar@{}[rd]|{(2)} & CH_n(N')_\mathbb{Q} \ar@{}[rd]|{(3)} \ar[d]^-{\wr} & CH_{n-c}(X')_\mathbb{Q} \ar[l]_-{\pi'^*} \ar[d]^-{\wr}\\
    H^{BM}_{2n,n}(Y') \ar[r]^-{sp} \ar[rd]_-{R_f(i)} & H^{BM}_{2n,n}(C') \ar[r]^-{i_*} \ar@{}[d]|{(4)} & H^{BM}_{2n,n}(N') & H^{BM}_{2n-2c,n-c}(X'). \ar[l]_-{\pi'^*}\\
      & H^{BM}_{2n-2c,n-c}(X') \ar[ru]^-{\pi'^*} \ar@{=}[rru] & &
  }
$$
We only need to apply the results already established: indeed the four squares $(1)$, $(2)$, $(3)$ and $(4)$ follow respectively from Proposition~\ref{compat_sp}, Proposition~\ref{compat_proper_pf}, Proposition~\ref{comp_smooth_pullback} and Corollary~\ref{refined_gysin_compat}.

\begin{proposition}
\label{tensor_prod}
Let $X$ and $Y$ be two separated $k$-schemes of finite type. Then the following diagram is commutative
$$
  \xymatrix{
  H^{BM}_{2m,m}(X)\otimes H^{BM}_{2n,n}(Y) \ar[r]^-{\otimes} \ar[d]^-{\wr} & H^{BM}_{2m+2n,m+n}(X\times Y) \ar[d]^-{\wr} \\
 CH_m(X)_\mathbb{Q}\otimes CH_n(Y)_\mathbb{Q} \ar[r]^-{\times} & CH_{m+n}(X\times Y)_\mathbb{Q},
  }
$$
where the map $\times$ in the lower row is the exterior product on Chow groups (\cite[1.10]{Ful}).
\end{proposition}
\proof

The proof uses Lemma~\ref{elem_BM} $(\ref{comp_prod})$ and is similar to Proposition~\ref{compat_proper_pf}.

\endproof

\subsection{Heart of weight structure on motives}\

Let $S$ be a quasi-projective $k$-scheme. Consider the two following categories:
\begin{itemize}
\item Let $CHM(S)$ be the category of Chow motives over the base $S$ with rational coefficients, defined in \cite[Definition 2.8]{CH}. This category is constructed from the category of projective $S$-schemes that are smooth over $k$. 
\item Let $Chow(S)$ be the smallest subcategory of $DM_{{\cyrrm{B}},c}(S)$, pseudo-abelian and stable by finite sums, that contains all elements of the form 
$$
p_!\mathbbold{1}_X(r)[2r]=M^{BM}(X/S)(r)[2r],
$$
where $X$ is a regular $k$-scheme 
endowed with a projective morphism $p:X\to S$.\footnote{H\'ebert showed that we can replace projective $S$-schemes by all proper $S$-schemes (\cite[Lemme 3.1]{Heb}).} 
Following \cite[Theorem 2.1.1]{Bon2} and \cite[Th\'eor\`eme 3.3]{Heb}, we know that $Chow(S)$ is the heart of a canonical weight structure on $DM_{{\cyrrm{B}},c}(S)$. 
\end{itemize}
Our aim is to identify the two categories above. There is a natural map
\begin{align*}
\mathcal{F}:CHM(S)&\to Chow(S)\\
(X,r)&\mapsto M^{BM}(X/S)(r)[2r],
\end{align*}
for any projective $S$-scheme $X$ that is smooth over $k$. 
The map $\mathcal{F}$ extends to the whole $CHM(S)$, since both categories are pseudo-abelian and stable by finite sums.
\begin{theorem}
\label{main_thm}
The map $\mathcal{F}$ is a functor, which is an equivalence between the two categories $CHM(S)$ and $Chow(S)$.
\end{theorem}
\proof
First we show that $\mathcal{F}$ is a functor. Let $X$ and $Y$ be projective $S$-schemes that are smooth over $k$. Using Lemma~\ref{Hom_BM} and Corollary~\ref{BM=Chow}, we have an isomorphism
\begin{align*}
\sigma_{X,Y}:Hom(M^{BM}(X/S)(r)[2r],M^{BM}(Y/S)(s)[2s])
&\overset{\epsilon_{X,Y}}{\simeq} H^{BM}_{2d_Y+2r-2s,d_Y+r-s}(X\times_SY)\\
&\simeq CH_{d_Y+r-s}(X\times_SY)_\mathbb{Q},
\end{align*}
which shows that $\mathcal{F}$ induces a bijection on morphisms.
It remains to see that $\mathcal{F}$ preserves composition of morphisms. On the one hand, Proposition~\ref{composition_law} shows that $\epsilon_{X,Y}$ gives the desired formal composition law on $CHM(S)$; on the other hand, it follows from Proposition~\ref{compat_proper_pf}, Proposition~\ref{compat_refined_gysin} and Proposition~\ref{tensor_prod} that the corresponding actions on Borel-Moore homology are compatible with the ones on Chow groups. Consequently, $\mathcal{F}$ is a functor which is fully faithful. Since every regular $k$-scheme is smooth, $\mathcal{F}$ is essentially surjective by the definition of $Chow(S)$, therefore a equivalence.
\endproof

\section{Lemmas}

\subsection{On Lemma~\ref{elem_BM}}\
\label{proof_elem_BM}

The aim of this section is to show Lemma~\ref{elem_BM}. In the whole section, $k$ is a perfect field.

\begin{definition}
\label{def_mot_coh}
Let $f:X\to Spec\ (k)$ be a smooth morphism.
The \textbf{(homological) motive} of $X$ over $k$ is defined as the object $M(X)=f_!f^!\mathbbold{1}_{k}$ in the category $DM_{{\cyrrm{B}},c}(k)$, and the \textbf{motivic cohomology} of $X$ is defined as
$$
H^{p,q}(X)\simeq Hom_{DM_{{\cyrrm{B}},c}(k)} (M(X),\mathbbold{1}_{k}(q)[p])
$$
\end{definition}


The purity isomorphism implies that motivic cohomology for smooth $k$-schemes agrees with Borel-Moore homology: 
\begin{lemma}
\label{ident_BM_mot}
For any $d$-dimensional smooth $k$-scheme $X$, there are isomorphisms
$$
M(X)\simeq M^{BM}(X)(d)[2d]\textrm{\ and\ \ } H^{BM}_{p,q}(X)\simeq H^{2d-p,d-q}(X).
$$
\end{lemma}
As a counterpart of the basic functorialities of Borel-Moore motives (Lemma~\ref{bm_basic}), the motives have covariant functoriality with respect to all morphisms, contravariant functoriality with respect to projective morphisms, and Gysin triangle (\cite[Definition 2.20]{Deg2}).
The corresponding version of Lemma~\ref{elem_BM} in terms of motivic cohomology follows from \cite[3.15]{Deg2}:
\begin{lemma}
\label{elem_coh_mot}
\begin{enumerate}
\item Let $p:X\to Y$ be a finite surjective morphism between irreducible smooth $k$-schemes. Then the map $p_*:H^{0,0}(X)\to H^{0,0}(Y)$ on motivic cohomology of degree $(0,0)$ induced by the Gysin morphism is the multiplication by the degree of field extension $[k(X):k(Y)]$ as a map $\mathbb{Q}\to\mathbb{Q}$.

\item Let $f:X\to Y$ be a morphism between irreducible smooth $k$-schemes. Then the map $f^*:H^{0,0}(Y)\to H^{0,0}(X)$ induced by functoriality is the identity map $\mathbb{Q}\to\mathbb{Q}$.

\item Let $i:Z\to X$ be a closed immersion of codimension $1$ between irreducible smooth $k$-schemes (i.e. $Z$ is a divisor of $X$). Denote by $x$ the generic point of $X$. Then there is a commutative diagram
$$
  \xymatrix{
    H^{1,1}(X-Z) \ar[rr]^-{\partial_{X,Z}} \ar[rd] & & H^{0,0}(Z)=\mathbb{Q},\\
       &  \hat{H}^{1,1}(x)=\kappa(x)^*_{\mathbb{Q}} \ar[ru]^-{ord_Z} &
  }
$$
where $ord_Z:\kappa(x)^*_\mathbb{Q}\to\mathbb{Q}$ is the order of vanishing along $Z$.

\item Let $X$ and $Y$ be irreducible smooth $k$-schemes of dimension respectively $d$ and $d'$. Then there is a canonical isomorphism $M(X)\otimes M(Y)\simeq M(X\times Y)$ such that the induced map $H^{0,0}(X)\otimes H^{0,0}(Y)\to H^{0,0}(X\times Y)$ sends the class $1\otimes 1$ to the class $\sum_I1$, where elements of $I$ are irreducible components of $X\times Y$.
\end{enumerate}

\end{lemma}
Therefore in order to prove Lemma~\ref{elem_BM} we need to check the compatibility of basic functorialities and the K\"unneth formula via Lemma~\ref{ident_BM_mot}. The K\"unneth formula follows from \cite[1.1.36]{CD}, the smooth functoriality follows from the compatibility between purity isomorphisms (Corollary~\ref{purity_indep}), and by Lemma~\ref{uniqueness_h}, it remains to check the proper functoriality for a projective morphism:
\begin{lemma}
\label{BM_mot_proj_push}
Let $X$ and $Y$ be smooth $k$-schemes, $d=\textrm{dim\ }X$ and $p:Y\to X$ be a projective morphism of relative dimension $c$. 
Then the following diagram is commutative:
$$
  \xymatrix{
   M^{BM}(X) \ar[r]^-{p^*}  \ar[d]^-{\wr} & M^{BM}(Y) \ar[d]^-{\wr} \\
     M(X)(-d)[-2d] \ar[r]^-{p^*} & M(Y)(-c-d)[-2c-2d].
  }
$$
\end{lemma}

\proof

A projective morphism is the composition of a closed immersion and the projection of a projective bundle, therefore it is sufficient to discuss these two cases.

Suppose that $i:Z\to X$ is a closed immersion of codimension $c$ between regular schemes. Denote by $f:X\to S$ the structure morphism. According to \cite[Definition 2.20]{Deg2} the Gysin morphism $M(X)\to M(Z)(c)[2c]$ is defined by the composition $M(X)\to M_Z(X)\simeq M(Z)(c)[2c]$, where $M_Z(X)$ is canonically identified with $f_!i_*i^*f^!\mathbbold{1}_S$, the first map is induced by the adjunction $1\xrightarrow{ad(i^*,i_*)}i_*i^*$ and the second map is the purity isomorphism. Therefore we have the following diagram:
$$
  \xymatrix{
    f_!f^*\mathbbold{1}_S \ar[r]^-{ad(i^*,i_*)} \ar[d]^-{\wr}_{P(f)^{-1}} & f_!i_*i^*f^*\mathbbold{1}_S \ar[d]^-{\wr}_{P(f)^{-1}} \ar[rd]^-{P(f\circ i)^{-1}}_{\sim} & \\
     f_!f^!\mathbbold{1}_S(-d_X)[-2d_X] \ar[r]^-{ad(i^*,i_*)} & f_!i_*i^*f^!\mathbbold{1}_S(-d_X)[-2d_X] \ar[r]^-{P(i)^{-1}}_{\sim} & f_!i_!i^!f^!\mathbbold{1}_S(-d_Z)[-2d_Z],
  }
$$
which is commutative: the square on the left is clearly commutative, and the triangle on the right follows from Lemma~\ref{compose_purity}.

The projective bundle case is more complicated, and we give here an idea of the proof. By localization sequence, we are reduced to the case of a trivial projective bundle  $Y=\mathbb{P}^n_X$. By K\"unneth formula, we may suppose that $X=k$ is a point. Then we proceed by induction, as in \cite[Theorem 3.2]{Deg6}. Note that the case $n=1$ follows from the fact that the canonical proper functoriality map $\mathbbold{1}_k\to M^{BM}(\mathbb{P}^1)$ is a section of the first Chern class of the canonical bundle $c_1(\mathcal{O}(1)):M^{BM}(\mathbb{P}^1)\to \mathbbold{1}_k$.
\endproof

\subsection{Complements on purity isomorphisms}\

\label{proof_BC_proper_smooth}

In this section we work on some additional properties of purity isomorphisms in order to prove some lemmas in Chapter~\ref{chapter2}. Concerning the relative purity isomorphism, we recall the following properties from \cite{CD}:
\begin{recall}
\label{recall_DM2}
\begin{enumerate}
\item (Definition 1.1.2. and 2.4.20.) For any smooth morphism $f:Y\to X$ of relative dimension $d$, there is a functor $f_\#:DM_{{\cyrrm{B}},c}(Y)\to DM_{{\cyrrm{B}},c}(X)$ which is a left adjoint of the functor $f^*$, such that if $f$ is an open immersion, then $f_\#=f_!$; there is also a functor $\Sigma_f:DM_{{\cyrrm{B}},c}(Y)\to DM_{{\cyrrm{B}},c}(Y)$ depending only on $f$, with a canonical isomorphism $f_\#\simeq f_!\Sigma_f$;
\item (Corollary 2.4.37. and Definition 2.4.38.)
\label{sigma_compatible}
For any $K\in DM_{{\cyrrm{B}},c}(Y)$ there is an isomorphism $\Sigma_f(K)\simeq K(d)[2d]$, and
for any cartesian square
$$
  \xymatrix{
    Y' \ar[r]^-{q} \ar[d]_-{g} & X' \ar[d]^-{f} \\ 
    Y \ar[r]^-{p} & X,
  }
$$
where $f$ and $g$ are smooth, there is an isomorphism of functors $\Sigma_gq^*\simeq q^*\Sigma_f$
such that the following diagram of functors is commutative:
$$
  \xymatrix{
    \Sigma_gq^* \ar[r]^-{\sim} \ar[d]_-{\wr} & q^*(d)[2d]. \\
    q^*\Sigma_f \ar[ru]_-{\sim} & 
  }
$$
\item (Theorem 2.4.50.) By definition, the relative purity isomorphism is deduced from the composed isomorphism $f_\#\simeq f_!\Sigma_f\simeq f_!(d)[2d]$ by adjunction.
\end{enumerate}
\end{recall}

\noindent
\textit{Proof of Lemma~\ref{BC_proper_smooth}.}
We need to show that the following diagram is commutative
$$
  \xymatrix{
    f_!f^* \ar[r]^-{ad_{(p^*,p_*)}} \ar[d]_-{ad_{(q^*,q_*)}} \ar@{}[rd]|{(1)} & p_!p^*f_!f^* \ar[r]^-{\sim} \ar[d]^-{\wr} \ar@{}[rd]|{(2)} & p_!p^*f_\#f^*(-d)[-2d] \ar[r]^-{ad'_{(f_\#,f^*)}} \ar[d]^-{\wr} \ar@{}[rd]|{(3)} & p_!p^*(-d)[-2d] \\
    f_!q_!q^*f^* \ar[r]^-{\sim} & p_!g_!q^*f^*\ar[r]^-{\sim} & p_!g_\#q^*f^*(-d)[-2d] \ar@{=}[r] 
&  p_!g_\#g^*p^*(-d)[-2d]. \ar[u]_-{ad'_{(g_\#,g^*)}}
  }
$$
The squares $(1)$ and $(3)$ 
follow from a formal compatibility between adjunction pairs.
For the square $(2)$, divide the diagram as follows:
$$
  \xymatrix{
    g_\#q^* \ar[r]^-{\sim} \ar[d]_-{\wr}  & g_!\Sigma_gq^* \ar[r]^-{\sim} \ar[d]^-{\wr} \ar@{}[rd]|{(5)} & g_!q^*(d)[2d] \ar@{=}[d] \\
    p^*f_\# \ar[rd] \ar@{}[r]|{(4)} & g_!q^*\Sigma_f \ar[r]^-{\sim} \ar[d]^-{\wr} \ar@{}[rd]|{(6)} & g_!q^*(d)[2d] \ar[d]^-{\wr} \\
 & p^*f_!\Sigma_f \ar[r]^-{\sim} & p^*f_!(d)[2d].
  }
$$
The diagrams $(4)$ and $(6)$ are commutative by formal compatibility, and the square $(5)$ comes from Recall~\ref{recall_DM2} (\ref{sigma_compatible}).
\endproof

\noindent
\textit{Proof of Lemma~\ref{lemma_monoidal_BC}.}
The first assertion is a formal compatibility, and by Lemma~\ref{uniqueness_h} the third assertion follows from the two previous ones. Therefore it is sufficient to show the second assertion. Divide the diagram as
$$
  \xymatrix{
    q_!\mathbbold{1}_{X'}(d)[2d] \ar[d]^-{\wr} \ar[r]^-{\sim} & q_\#\mathbbold{1}_{X'}   \ar[r] \ar[d]^-{\wr} & \mathbbold{1}_{Y'} \ar@{=}[d] \\
   f^*p_!\mathbbold{1}_X(d)[2d] \ar[r]^-{\sim} & f^*p_\#\mathbbold{1}_X \ar[r] & f^*\mathbbold{1}_{Y}.
  }
$$
The square on the left follows from the square $(2)$ in the proof of Lemma~\ref{BC_proper_smooth}, and the square on the right follows from adjunction.

\endproof

Now consider the absolute purity case. Let $i:Z\to X$ be a closed immersion of codimension $c$ between regular schemes. Following \cite[Remark 2.3.3]{Deg5}, there is a map $i_*\mathbbold{1}_{Z}\to\mathbbold{1}_{X}(c)[2c]$, which by adjunction induces the absolute purity isomorphism $\mathbbold{1}_{Z}\simeq i^!\mathbbold{1}_{X}(c)[2c]$.
The following lemma follows from \cite[Corollary 2.4.4]{Deg5}:
\begin{lemma}
\label{compat_abs_purity}
Consider a cartesian square of regular schemes
$$
  \xymatrix{
    Z' \ar[r]^-{i'}  \ar[d]_-{g} & X' \ar[d]^-{f} \\
    Z\ar[r]^-{i} & X,
  }
$$
where $i$ and $i'$ are closed immersions of codimension $c$ and $f$ is transversal to $i$. Then the following diagram is commutative:
$$
  \xymatrix{
    f^*i_*\mathbbold{1}_Z \ar[r]^-{} \ar[d]^-{\wr} & f^*\mathbbold{1}_X(c)[2c] \ar@{=}[d] \\
    i'_*\mathbbold{1}_{Z'} \ar[r]^-{} & \mathbbold{1}_{X'}(c)[2c].
  }
$$
\end{lemma}

\noindent
\textit{Proof of Lemma~\ref{compat_purity}.}
First we show that there is a commutative diagram
$$
  \xymatrix{
    f_!i'_*g^*(d)[2d] \ar[r]^-{\sim} \ar[d]^{\wr} & f_\#i'_*g^* \ar[r]^-{\sim} \ar[d]^-{\wr} & f_\#f^*i_*  \ar[d]^-{ad'_{(f_\#,f^*)}} \\
    i_*g_!g^*(d)[2d] \ar[r]^-{\sim} & i_*g_\#g^* \ar[r]^-{ad'_{(g_\#,g^*)}} & i_*.
  }
$$
The square on the left 
is similar to the square $(2)$ in the proof of Lemma~\ref{BC_proper_smooth} (the intersection is transversal since $f$ is smooth), and the one on the right is a formal compatibility between adjunctions. On the other hand, using Lemma~\ref{compat_abs_purity}, we have the following big commutative diagram:
$$
  \xymatrix{
    f_!i'_*\mathbbold{1}_{Z'} \ar[r]^-{} \ar[d]_{\wr} & f_\#i'_*\mathbbold{1}_{Z'}(-d)[-2d] \ar[d]^-{\wr} \ar[r]^-{} \ar[rd]^-{\sim} & f_\#\mathbbold{1}_{X'}(c-d)[2c-2d] \ar[rd]^-{\sim} & \\
   i_*g_!\mathbbold{1}_{Z'} \ar[r]^-{} & i_*g_\#\mathbbold{1}_{Z'}(-d)[-2d] \ar[rd]^-{ad'_{(g_\#,g^*)}} & f_\#f^*i_*\mathbbold{1}_{Z}(-d)[-2d] \ar[d]^-{ad'_{(f_\#,f^*)}} \ar[r]^-{} & f_\#f^*\mathbbold{1}_{X}(c-d)[2c-2d] \ar[d]^-{ad'_{(f_\#,f^*)}} \\ 
 & & i_*\mathbbold{1}_Z(-d)[-2d] \ar[r]^-{} & \mathbbold{1}_X(c-d)[2c-2d],
  }
$$
which proves the result.
\endproof

\noindent
\textit{Proof of Lemma~\ref{compat_gysin_func}.}
By Lemma~\ref{uniqueness_h} it is sufficient to show the first assertion, and we may assume that $X=S$. From Lemma~\ref{compat_abs_purity} we know that the following diagram is commutative: 
$$
  \xymatrix{
    p^*i_*\mathbbold{1}_{Z} \ar[d]^-{\wr} \ar[r]^-{}& p^*\mathbbold{1}_X(c)[2c]  \ar@{=}[d] \\
   i'_*\mathbbold{1}_{Z'} \ar[r]^-{} & \mathbbold{1}_Z(c)[2c].
  }
$$
Applying the functor $i_*$ to the previous diagram, and we obtain the result by a formal computation of adjunctions.
\endproof

\end{document}